\pdfoutput=1
\documentclass[11pt,a4paper,headings=normal,numbers=noenddot,parskip=never,abstract=true]{scrartcl}
\usepackage[T1]{fontenc}
\usepackage[utf8]{inputenc}
\usepackage{lmodern,microtype}
\usepackage[a4paper,margin=1.22in,includehead=false,includefoot=false]{geometry}
\usepackage{amsmath,amsfonts,amssymb,amsthm,mathtools,amscd}
\usepackage{tikz-cd}
\usepackage{booktabs,array,framed,xcolor}
\usepackage[only,longmapsfrom]{stmaryrd}
\usepackage{scrlayer-scrpage}
\usepackage[colorlinks=true,linkcolor=blue!40!black,citecolor=blue!40!black,urlcolor=blue!40!black]{hyperref}
\setkomafont{disposition}{\normalfont\bfseries}
\setkomafont{title}{\normalfont\Large}
\setkomafont{section}{\normalfont\Large\bfseries}
\setkomafont{subsection}{\normalfont\large\bfseries}
\clearpairofpagestyles
\definecolor{reviewgray}{RGB}{245,245,244}
\definecolor{reviewblue}{RGB}{237,243,249}

\newtheorem{theorem}{Theorem}[section]
\newtheorem{proposition}[theorem]{Proposition}
\newtheorem{lemma}[theorem]{Lemma}
\newtheorem{corollary}[theorem]{Corollary}
\newtheorem{mainA}{Theorem}

\theoremstyle{definition}
\theoremstyle{remark}
\newcommand{\Q}{\mathbb Q}\newcommand{\Z}{\mathbb Z}
\newcommand{\R}{\mathbb R}\newcommand{\C}{\mathbb C}
\newcommand{\OP}{\mathbb OP}\newcommand{\HP}{\mathbb HP}
\newcommand{\Spin}{\mathrm{Spin}}\newcommand{\String}{\mathrm{String}}
\newcommand{\Sp}{\mathrm{Sp}}\newcommand{\Diff}{\mathrm{Diff}}
\newcommand{\Ahat}{\widehat A}\newcommand{\ch}{\operatorname{ch}}
\newcommand{\ph}{\operatorname{ph}}\newcommand{\Ind}{\operatorname{Ind}}
\newcommand{\hofib}{\operatorname{hofib}}\newcommand{\im}{\operatorname{im}}
\newcommand{\scal}{\operatorname{scal}}\newcommand{\Fr}{\operatorname{Fr}}
\newcommand{\cR}{\mathcal R}\newcommand{\cC}{\mathcal C}
\newcommand{\cP}{\mathcal P}\newcommand{\rst}{\mathrm{rst}}
\newcommand{\Hh}{\mathbb H}\newcommand{\Hom}{\operatorname{Hom}}
\newcommand{\op}{\mathrm{op}}
\newcommand{\Bpsc}{B\cP_{\Spin(16)}^{2,\rst}}
\newcommand{\Bord}{B\cC_{\Spin(16)}}
\newcommand{\Bdiff}{B\widehat{\Diff}_0(\OP^2)}
\DeclareMathOperator{\Sym}{Sym}
\DeclareMathOperator{\Wit}{Wit}
\hypersetup{pdftitle={Cayley plane bundles, fibrewise positive scalar curvature, and the Witten genus},
pdfsubject={Positive scalar curvature and family indices of Cayley plane bundles}}
\title{Cayley plane bundles, fibrewise positive scalar curvature, and the Witten genus\vspace{4pt}}
\author{Yuanchu Li}
\date{}

\begin{document}
\maketitle
\begingroup
\renewcommand{\thefootnote}{}
\footnotetext{Email address:
	\href{mailto:yuanchuli@u.nus.edu}{\nolinkurl{yuanchuli@u.nus.edu}}.}
\endgroup
\begin{abstract}

By constructing a smooth Cayley plane bundle over $S^8$
with fibrewise positive scalar curvature and nonzero Witten genus, 
we show that the rational kernel of
the Dirac index in String bordism is spanned by total
spaces of Cayley plane bundles
with smooth fibrewise positive scalar curvature.

\end{abstract}

\tableofcontents

\section{Introduction}\label{sec:introduction}

In this paper, we construct a smooth Cayley plane bundle over $S^8$ with fibrewise positive scalar curvature whose String total space has nonzero Witten genus. In particular, the bundle has no reduction to the standard $F_4$ action.
Together with homogeneous Cayley bundles, this construction gives geometric representatives for the rational kernel of the ordinary Dirac index in String bordism.

\subsection*{Cayley bundles as representatives in String bordism}

The Cayley plane is the homogeneous space
$\OP^2=F_4/\Spin(9)$, equipped with its standard
homogeneous metric $g_{\mathrm{std}}$ of positive
sectional curvature.
We call a smooth Cayley bundle \emph{homogeneous}
if its structure group reduces to $\operatorname{Isom}(\OP^2,g_{\mathrm{std}})
=F_4$. The universal homogeneous Cayley bundle is then given by
\begin{equation}\label{eq:homogeneous-bundle}
	\OP^2\longrightarrow B\Spin(9)\longrightarrow BF_4.
\end{equation}

The Witten genus is a ring homomorphism
from the String bordism to the graded ring
of modular forms, as we recall in Section~\ref{subsec:indices}. 
Using Cayley bundles, Jung and Dessai identified
the rational kernel of the Witten genus as
\begin{equation}\label{eq:homogeneouskernel}
	\ker\big(
	\Omega_*^{\String}
	\xrightarrow{\Wit}
	MF_*
	\big)_\Q
	=
	\operatorname{span}_{\Q}
	\left\{
	\begin{array}{l}
		\text{String bordism classes of total spaces}\\
		\text{of \textit{homogeneous Cayley bundles}}
	\end{array}
	\right\},
\end{equation}
and McTague extended this result to
$\mathbb Z[1/6]$, see \cite{Des09} and \cite{McT14}. 
This geometric approach continues a line of work of Stolz and
Kreck--Stolz \cite{St92,KS93}, where they characterized
the kernel of the real Dirac index
$\alpha^{\Spin}_{KO}$ on Spin bordism
in terms of total spaces of $\HP^2$-bundles
associated to the standard homogeneous
$\operatorname{PSp}(3)$-action.

\medskip

Let $\cR^+_{\mathrm{std}}(\OP^2)$ denote the path component of
$g_{\mathrm{std}}$ in the space of psc metrics.
Using the Spin diffeomorphism group
$\widehat{\Diff}_0(\OP^2)$
of Lemma~\ref{lem:spin}, we define the homotopy quotient 
\[
\mathcal M^+_{\mathrm{std}}
=
E\widehat{\Diff}_0(\OP^2)
\times_{\widehat{\Diff}_0(\OP^2)}
\cR^+_{\mathrm{std}}(\OP^2),
\]
which classifies smooth Cayley bundles
with structure group $\Diff_0(\OP^2)$, the specified
vertical Spin structure, and fibrewise psc metrics
in the standard component.

Since the standard $F_4$ action lifts to
$\widehat{\Diff}_0(\OP^2)$ and fixes
$g_{\mathrm{std}}$, every homogeneous Cayley bundle carries
an induced vertical Spin structure and
a fibrewise psc metric isometric to
$g_{\mathrm{std}}$ on each fibre.
Consequently, we have maps of classifying spaces
\begin{equation}\label{eq:homogeneous-psc}
	BF_4
	\longrightarrow
	\mathcal M^+_{\mathrm{std}}
	\longrightarrow
	\Bdiff.
\end{equation}

In this paper, we establish a parallel geometric description
for the larger rational kernel of the ordinary Dirac orientation,
passing from \textit{homogeneous} families classified by $BF_4$ to \textit{fibrewise psc} families
classified by $\mathcal M^+_{\mathrm{std}}$. More precisely, we prove
\begin{equation}\label{eq:psc-kernel}
	\ker\big(
	\Omega_*^{\String}
	\xrightarrow{\alpha_{KO}^{\String}}
	KO_*
	\big)_\Q
	=
	\operatorname{span}_{\Q}
	\left\{
	\begin{array}{l}
		\text{String bordism classes of total spaces}\\
		\text{of \textit{fibrewise psc Cayley bundles}}
	\end{array}
	\right\}.
\end{equation}
Here $\alpha_{KO}^{\String}$ denotes the ordinary real
Dirac index on String bordism, whose complexification
is given by the $\widehat A$-genus, the constant term
of the Witten genus. In fact, this equality also holds if the fibrewise psc
metrics are required to lie in the standard path component.

\subsection*{Main results}

For a Cayley bundle $\pi:E\to X$ with String total space,
the $7$-connectedness of the fibre allows the
underlying Spin structure on $E$ to be expressed as
the sum of the natural vertical Spin structure on $T_\pi E$ (see Lemma~\ref{lem:spin})  and a Spin structure on $X$. The ordinary Dirac index of $E$ is then given by the functoriality of $KO$-pushforwards
\begin{equation}\label{eq:psc-total-index}
	\alpha_{KO}^\String(E)
	=(E\to *)_!^{KO}(1)=(X\to *)_!^{KO}\circ\pi_!^{KO}(1).
\end{equation}
By the Atiyah--Singer family index theorem \cite{AS71V,Hit74}, we have
\begin{equation}\label{eq:family-index}
	\Ind_{KO}(D_\pi)=\pi_!^{KO}(1),
\end{equation}
where $\Ind_{KO}(D_\pi)$ denotes the $KO$-valued analytic family index of the vertical Dirac operators $D_\pi$.
When $E$ admits fibrewise psc metrics, $D_\pi$ are invertible by the Lichnerowicz formula, hence 
$\Ind_{KO}(D_\pi)$ vanishes. 
Combining \eqref{eq:psc-total-index} and
\eqref{eq:family-index}, we have established the
inclusion from right to left in \eqref{eq:psc-kernel} 
for bundles with structure group $\Diff_0(\OP^2)$.

For the reverse inclusion, our main construction enters as follows. 

\begin{mainA}\label{thm:A}
	There is a smooth Cayley plane bundle $E_+$ over $S^8$ with a smooth fibrewise psc metric in
	$\cR^+_{\mathrm{std}}(\OP^2)$, such that its Witten genus is a nonzero multiple of the discriminant form $\Delta$. In particular, $E_+$ is not homogeneous.
\end{mainA}

Theorem~\ref{thm:A} also illustrates how twisted family
indices can detect nonhomogeneity in the presence of
fibrewise positive scalar curvature.
For related work on fibrewise curvature and the topology
of smooth bundles, see \cite{BHSW10,MS26}.
For characteristic numbers arising from twisted Dirac
indices on String manifolds, see also \cite{HH22}.

\medskip

Returning to \eqref{eq:psc-kernel}, let $M$ be a closed
String manifold with torsion ordinary Dirac index. Its Witten genus $\Wit(M)$ is a cusp form and hence divisible by $\Delta$ in the ring of modular forms. 
Since $\Wit(E_+)$ is a nonzero
multiple of $\Delta$, we may write $\Wit(M)=\Wit(E_+)f$ for a
modular form $f$ with rational coefficients. By rational
surjectivity of the Witten genus \cite{McT14}, there exist a
positive integer $k$ and a closed String manifold $Y$ such that
$\Wit(Y)=kf$. It follows from the multiplicativity of the Witten genus that
\[
\Wit(kM)=\Wit(E_+\times Y).
\]
By \eqref{eq:homogeneouskernel}, the class
$k[M]-[E_+\times Y]$ is a rational linear combination of String
bordism classes of total spaces of homogeneous Cayley bundles.
These bundles carry their standard fibrewise psc metrics,
while the product bundle
\[
E_+\times Y\longrightarrow S^8\times Y
\]
carries the pullback fibrewise psc metric in the standard component. Thus $[M]$ lies in
the rational span of String bordism classes of total spaces
of fibrewise psc Cayley bundles.
Together with the opposite inclusion, this proves the
following corollary.
\begin{corollary}\label{cor:kernel}
	The rational kernel of the ordinary Dirac index on
	String bordism admits the following geometric description
	\begin{equation}
		\ker\big(
		\Omega_*^{\String}
		\xrightarrow{\alpha_{KO}^{\String}}
		KO_*
		\big)_\Q
		=
		\operatorname{span}_{\Q}
		\left\{
		\begin{array}{l}
			\text{String bordism classes of total spaces}\\
			\text{of \emph{fibrewise psc Cayley bundles}}
		\end{array}
		\right\}.
	\end{equation}
	In particular, the statement is unchanged when one restricts to bundles classified by
	$\mathcal M^+_{\mathrm{std}}$.
\end{corollary}

\subsection*{Outline of the construction}

The following sections are devoted to the construction in
Theorem~\ref{thm:A}. In Section~\ref{sec:objects}, we
apply Frenck's results \cite{Fre25} to construct a smooth Cayley
family represented by $\xi$, with vanishing ordinary family
index and nonzero first twisted family index.
It then remains to equip fibrewise psc metrics to some positive multiple of $\xi$.

In Section~\ref{sec:metrics}, we first express the normalized
family Pontryagin--Thom class $z_0$ of $\xi$ rationally as
the image of a class $y$ under a fixed-rank geometric
transfer. Using the homological stability theorem of
Galatius and Randal-Williams \cite{GRW17}, we realize a
positive multiple of the rational Hurewicz image of $y$
by a smooth $\theta$-structured family over a
closed oriented eight-manifold.
The disk bundle caps associated to this $\theta$-family give a family $j$ of psc loops in the psc cobordism category, and in turn produce a rational psc homotopy class $a_\Q$, whose forgetful image coincides with $z_0\otimes 1$ under
the Pontryagin--Thom equivalence.

In Section~\ref{sec:based}, we consider the punctured
bundle $\xi_\partial$ and compare its Pontryagin--Thom
class with $a_\Q$.
By applying the comparison square of Ebert and
Randal-Williams \cite{ERW22}, we obtain a continuous section of the
metric bundle of some multiple of the punctured family.
Capping and smoothing this section produces the desired
metrics in Theorem~\ref{thm:A}.

\section{Smooth Cayley plane bundles}\label{sec:objects}

In this section, we begin by discussing the vertical Spin structures
and characteristic classes of smooth Cayley plane bundles. 
Then we construct a smooth Cayley plane bundle over $S^8$ with trivial ordinary family index and nonzero Witten genus.

\subsection{Spin structures and family indices}\label{subsec:indices}
The Cayley plane $\OP^2$ is $7$-connected and has integral
cohomology ring
\[
H^*(\OP^2;\Z)\cong\Z[u]/(u^3),\qquad |u|=8,
\]
and hence admits unique Spin and String structures.
With the standard orientation and a suitable choice of $u$ as in \cite{BH58}, we have 
\begin{equation}\label{eq:fibre}
 \langle u^2,[\OP^2]\rangle=1,\qquad
 p_2(T\OP^2)=6u,\quad p_4(T\OP^2)=39u^2,\quad e(T\OP^2)=3u^2.
\end{equation}
In particular, $\operatorname{sig}(\OP^2)=1$ and $\Ahat(\OP^2)=0$.

\medskip

Let $\Diff_0(\OP^2)$ denote the identity component of the
diffeomorphism group of $\OP^2$.
The Spin lifts of its elements form a double cover of
$\Diff_0(\OP^2)$. In fact, we have that

\begin{lemma}\label{lem:spin}
	The Spin double cover of $\Diff_0(\OP^2)$
	is trivial, with its identity component denoted by $\widehat{\Diff}_0(\OP^2)$. Moreover, the homeomorphism
	\[
	\Diff_0(\OP^2)
	\xrightarrow{\;\cong\;}
	\widehat{\Diff}_0(\OP^2)
	\]
	induces natural vertical Spin structures on bundles classified by $B\Diff_0(\OP^2)$.  
\end{lemma}

\begin{proof}
	Let $\pi:E_\gamma\to S^2$ be the bundle obtained
	by clutching along an arbitrary loop $\gamma$ in
	$\Diff_0(\OP^2)$ based at the identity.
	The Serre spectral sequence gives
	$w_2(T_\pi E_\gamma)=\pi^*\omega$ for some
	$\omega\in H^2(S^2;\Z/2)$ and a class
	$u_\gamma\in H^8(E_\gamma;\Z)$ restricting to $u$.
	The Spin monodromy of the clutching loop is then computed by
	$\langle\omega,[S^2]\rangle$.
	The mod $2$ reduction $\bar u_\gamma$ satisfies
	$\pi_!(\bar u_\gamma^2)=1$ and $\mathrm{Sq}^1\bar u_\gamma=0$.
	Since $w_2(TS^2)=0$, the Wu and Cartan formulas give
	\[
	\begin{aligned}
		\langle\omega,[S^2]\rangle
		&=\langle w_2(TE_\gamma)\bar u_\gamma^2,[E_\gamma]\rangle\\
		&=\langle\mathrm{Sq}^2(\bar u_\gamma^2),[E_\gamma]\rangle\\
		&=\langle(\mathrm{Sq}^1\bar u_\gamma)^2,[E_\gamma]\rangle=0.
	\end{aligned}
	\]
	Thus every loop has trivial Spin monodromy, so the
	double cover is trivial. 
\end{proof}

For a real vector bundle $V$ of rank $d$, define the reduced Witten
bundle \cite[Section~0.2]{HHLZ21} by the formal $q$-series
\[
\Theta_q(V_\C)=\bigotimes_{r\ge1}\Sym_{q^r}(V_\C-\C^d).
\]
For a closed oriented manifold $M$, the Witten genus is
\[
\Wit(M)=\langle\Ahat(TM)\ch\Theta_q(T_\C M),[M]\rangle.
\]
If $M$ is Spin, this is the complex Dirac
index twisted by $\Theta_q(T_\C M)$.
If $M$ is String, $\Wit(M)$ is a modular form of weight
$\dim M/2$. 
Over $\Q$, the ideal of cusp forms is generated by
the discriminant form $\Delta$.

Let $\pi:E\to\Bdiff$ be the universal smooth Cayley bundle
with its specified vertical Spin structure, and let $D_\pi$
be the family of Spin Dirac operators along the fibres.
Using Bott periodicity,
we define
\begin{equation}\label{eq:indices}
	J(q)=\Ind_{KU}(D_\pi\otimes\Theta_q(T_{\pi,\C}))
	=\sum_{n\ge0}J_nq^n\in KU^0(\Bdiff)[[q]].
\end{equation}
Since $\Theta_q(T_{\pi,\C})=1+(T_{\pi,\C}-\C^{16})q+O(q^2)$,
we have
\[
\begin{aligned}
	J_0&=\Ind_{KU}(D_\pi)=\bigl(\Ind_{KO}(D_\pi)\bigr)_\C,\\
	J_1&=\Ind_{KU}(D_\pi\otimes(T_{\pi,\C}-\C^{16})).
\end{aligned}
\]
Applying the family index theorem \cite{AS71V}, we obtain by naturality
\begin{equation}\label{eq:familyindex}
	\ch J(q)=\pi_!\bigl(\Ahat(T_\pi)\ch\Theta_q(T_{\pi,\C})\bigr).
\end{equation}
Here $\Ind_{KO}$ and $\Ind_{KU}$ denote the real and
complex analytic family indices, respectively.
Define
\begin{equation}\label{eq:C}
	C=\ch_{[8]}J_1\in H^8(\Bdiff;\Q),
\end{equation}
where $\ch_{[r]}$ denotes the cohomological degree-$r$
component of the Chern character.

Let $\pi:E\to X$ be a smooth Cayley bundle classified by
a map $X\to\Bdiff$, where $X$ is a closed Spin manifold.
By the splitting $TE\cong T_\pi E\oplus\pi^*TX$,
multiplicativity of $\Ahat$ and $\Theta_q$,
\eqref{eq:familyindex}, and the projection formula, we obtain
\begin{equation}\label{eq:total}
	\Wit(E)=
	\langle
	\Ahat(TX)\ch\Theta_q(T_\C X)\ch J(q),[X]
	\rangle,
\end{equation}
where $J(q)$ and $C$ also denote their pullbacks along
the classifying map of the bundle.
In particular, since $TS^8$ is stably trivial, a Cayley
bundle over $S^8$ satisfies
\[
[q]\Wit(E)=\langle C,[S^8]\rangle.
\]

\medskip 

To relate the family indices to the Pontryagin--Thom construction,
we recall the definition of generalized Miller--Morita--Mumford
classes.
Let $\gamma_{16}\to B\Spin(16)$ be the universal Spin bundle
and $MT\Spin(16)=\mathrm{Th}(-\gamma_{16})$ the
Madsen--Tillmann spectrum.
For a smooth bundle $\rho:E\to X$ of closed $16$-manifolds
with a vertical Spin structure, the family
Pontryagin--Thom map
\[
	\alpha_\rho:X\longrightarrow\Omega^\infty MT\Spin(16)
\]
is adjoint to the composite
\[
\Sigma^\infty X_+
\longrightarrow\mathrm{Th}(-T_\rho E)
\longrightarrow MT\Spin(16),
\]
given by the fibrewise collapse and the classifying map
of the vertical Spin bundle.
For $t\in H^{16+s}(B\Spin(16);\Q)$ with $s>0$,
the universal generalized Miller--Morita--Mumford
(MMM) class $\kappa_t$ is the image
of $t$ under
\begin{equation}\label{eq:universal-mmm}
	H^{16+s}(B\Spin(16);\Q)
	\xrightarrow{\mathrm{Thom}}H^s(MT\Spin(16);\Q)
	\xrightarrow{\sigma}H^s(\Omega^\infty MT\Spin(16);\Q),
\end{equation}
where $\sigma$ is the cohomology suspension.
By the construction of the Gysin map, we have
\cite[Section~1.1]{GRW14MMM}
\begin{equation}\label{eq:family-mmm}
	\alpha_\rho^*\kappa_t
	=\rho_!t(T_\rho E)\in H^s(X;\Q).
\end{equation}
As cohomology suspensions, the positive degree
universal MMM classes are primitive, hence
invariant under translation in the infinite loop space.

For the universal Cayley bundle, denote its family
Pontryagin--Thom map by
\begin{equation}\label{eq:family-PT}
	\alpha_{\OP^2}:\Bdiff\longrightarrow
	\Omega^\infty MT\Spin(16).
\end{equation}
By \eqref{eq:familyindex} and \eqref{eq:family-mmm},
the class $C$ has the expression
\[
C=\alpha_{\OP^2}^*
\kappa_{
	[\Ahat(\gamma_{16})
	\ch(\gamma_{16,\C}-\C^{16})]_{24}}.
\]

\subsection{Construction of the smooth family}\label{sec:smooth}

This subsection is devoted to the following construction: 
\begin{proposition}\label{prop:smoothsphere}
	There is a smooth $\OP^2$-bundle over $S^8$ classified by
	\[
	\xi\in\pi_8\Bdiff
	\]
	with zero ordinary family index and nonzero pairing
	with $C$. 
	In addition, this bundle admits a
	fibrewise homotopy equivalence to the trivial $\OP^2$-bundle. 
\end{proposition}

\begin{proof}
	Let
	\[
	\mathcal S_{\mathrm{sm}}
	=\hofib\bigl(B\Diff^+(\OP^2)\longrightarrow
	B\operatorname{hAut}^+(\OP^2)\bigr)
	\]
	be the homotopy fibre at the trivial bundle. It classifies smooth
	oriented $\OP^2$-bundles equipped with a fibrewise homotopy
	equivalence to the trivial bundle.

	Following \cite{Fre25}, consider the subspace
	\begin{equation}\label{eq:smooth-bordism-image}
		\operatorname{Fib}^h_{\OP^2,8}
		=\im\bigl(
		\pi_8\mathcal S_{\mathrm{sm}}\otimes\Q
		\longrightarrow\Omega_{24}^{SO}\otimes\Q
		\bigr),
	\end{equation}
	spanned by the bordism classes of total spaces of
	fibre homotopy trivial Cayley plane bundles over $S^8$.
	By \eqref{eq:fibre}, all Pontryagin classes of $\OP^2$ are
	generated by $p_2$, with $p_2^2\ne0$ and $p_2^3=0$.
	Moreover, $\OP^2$ is $7$-connected and $16$-dimensional, so by \cite[Remark~1.1(i)(b)]{Fre25} the degree $8$ lies in the unblocking range $13=\min(16-1,2\cdot7-1)$, and
	\cite[Corollary~C(ii)]{Fre25} shows that
	\[
	\dim\operatorname{Fib}^h_{\mathbb{OP}^2,8}=2.
	\]

	For such a bundle $\pi:E\to S^8$, fix a basepoint $*\in S^8$
	and an identification $E_*\cong\OP^2$, and let $i:\OP^2\cong E_*\hookrightarrow E$
	be the fibre inclusion.
	Choose a fibrewise homotopy equivalence such that the composition
	\[
	r: E\xrightarrow{\simeq}
	S^8\times\OP^2
	\xrightarrow{\operatorname{pr}_2}
	\OP^2
	\]
	satisfies $r\circ i\simeq\mathrm{id}_{\OP^2}$.
	It identifies
	\[
	H^*(E;\Q)\cong\Q[t,u]/(t^2,u^3),
	\qquad |t|=|u|=8,\qquad
	\langle tu^2,[E]\rangle=1,
	\]
	where $t$ is the generator of
	$H^8(S^8;\Q)$ that evaluates to $1$, and
	$u$ is obtained from
	\eqref{eq:fibre}. 
	We then have 
	\[
	\ker i^*=(t)=\Q\{t,tu,tu^2\}.
	\]
	Since $i^*T_\pi E\cong T\OP^2$ and
	$r\circ i\simeq\mathrm{id}_{\OP^2}$, we have $\ph(T_\pi E)-r^*\ph(T\OP^2)\in \ker i^*$.
	Thus there are unique $a,b,c\in\Q$ such that
	\begin{equation}\label{eq:tangentdefect}
		\ph(T_\pi E)-r^*\ph(T\OP^2)
		=at+btu+ctu^2.
	\end{equation}
	Since $TS^8$ is stably trivial, we have
	$p_i(TE)=p_i(T_\pi E)$, which we simply denote by $p_i$.
	Newton identities and \eqref{eq:fibre} give
	\[
	\begin{aligned}
		\ph(T_\pi E)
		&=16-\frac{p_2}{6}
		+\frac{p_2^2-2p_4}{10080}
		+\frac{-p_2^3+3p_2p_4-3p_6}{119750400},\\
		r^*\ph(T\OP^2)
		&=16-u-\frac{u^2}{240}.
	\end{aligned}
	\]
	Comparing with \eqref{eq:tangentdefect},
	we obtain
	\[
	\begin{aligned}
		p_2&=6u-6at,\\
		p_4&=39u^2-(36a+5040b)tu,\\
		p_6&=(-234a-30240b-39916800c)tu^2.
	\end{aligned}
	\]
	Then we compute the only possibly nonzero Pontryagin numbers of $E$ as follows
	\[
	\begin{aligned}
		\langle p_2^3,[E]\rangle&=-648a,\\
		\langle p_2p_4,[E]\rangle&=-450a-30240b,\\
		\langle p_6,[E]\rangle&=-234a-30240b-39916800c.
	\end{aligned}
	\]
	The three Pontryagin numbers above determine the rational oriented
	bordism class of $E$.
	Moreover, they are related to $(a,b,c)$ by an
	invertible matrix.
	Thus $\operatorname{Fib}^h_{\mathbb{OP}^2,8}$
	is identified with a two-dimensional subspace of
	$\Q^3$ in the coordinates $(a,b,c)$.

	We now choose a class in this two-dimensional subspace with vanishing ordinary index but
	nonzero first twisted index.
	Consider the characteristic numbers
	\[
	\begin{aligned}
		\ell_L&=\langle L_{24}(TE),[E]\rangle,\\
		\ell_{\Ahat}&=\langle\Ahat_{24}(TE),[E]\rangle,\\
		\ell_1&=\big\langle
		[\Ahat(TE)\ch(T_\C E-\C^{24})]_{24},[E]\big\rangle,
	\end{aligned}
	\]
	which may be regarded as linear functionals on $\Q^3$ defined in terms of $(a,b,c)$.

	Since every bundle representative in $\operatorname{Fib}^h_{\mathbb{OP}^2,8}$ has total space
	homotopy equivalent to $S^8\times\OP^2$, its signature vanishes. 
	Since $\ell_L(1,0,0)=-14/15\ne0$,
	$\ker\ell_L$ has dimension two.
	It follows that, under the coordinates $(a,b,c)$,
	\[
	\operatorname{Fib}^h_{\OP^2,8}=\ker\ell_L.
	\]
	A direct calculation with the characteristic series further shows $\ell_{\Ahat}$ and
	$\ell_1$ depend only on $(b,c)$, and linearly independent on $\ker\ell_L$.
	Thus we may choose a class in $\operatorname{Fib}^h_{\mathbb{OP}^2,8}$ with
	\begin{equation}\label{eq:functional-condition}
		\ell_{\Ahat}=0, \qquad \ell_1\ne0.
	\end{equation}

	Taking a preimage under \eqref{eq:smooth-bordism-image}
	in $\pi_8\mathcal S_{\mathrm{sm}}\otimes\Q$ and clearing denominators, we obtain
	an actual smooth bundle $\pi:E\to S^8$ with a
	fibrewise homotopy equivalence to the trivial bundle,
	still satisfying \eqref{eq:functional-condition}. 
	Since $S^8$ is simply connected, its structure group reduces
	to $\Diff_0(\OP^2)$, and Lemma~\ref{lem:spin} gives its
	specified vertical Spin lift. Denote the resulting class by
	\[
	\xi\in \pi_8\Bdiff.
	\]

	By \eqref{eq:familyindex} and the stable triviality of
	$TS^8$, the complexification of $\Ind_{KO}(D_\pi)$ has
	rank $\Ahat(\OP^2)=0$ and degree-eight Chern-character
	pairing $\ell_{\Ahat}=0$.
	Since complexification is rationally injective and
	$KO^{-16}(S^8)$ is torsion-free,
	we conclude that $\Ind_{KO}(D_\pi)=0$.
	The same family index formula also identifies
	$\langle C,\xi\rangle$ with the nonzero number $\ell_1$.
\end{proof}

\section{From the Pontryagin--Thom class to a rational psc lift}\label{sec:metrics}\label{sec:transfer}\label{sec:finite}

As in \eqref{eq:family-PT}, let
\[
\alpha_{\OP^2}:\Bdiff
\longrightarrow\Omega^\infty MT\Spin(16)
\]
be the family Pontryagin--Thom map for $\OP^2$-bundles
with their specified vertical Spin structures.
Subtracting the value of $\alpha_{\OP^2}$ at the trivial bundle
gives the normalized class
\begin{equation}\label{eq:z}
	z_0:=
	\bigl(\alpha_{\OP^2}(-)-\alpha_{\OP^2}(\OP^2)\bigr)_*\xi
	\in\pi_8MT\Spin(16).
\end{equation}
We pair characteristic classes with homotopy classes
of a Thom spectrum using the stable Hurewicz map
and the Thom isomorphism.
By \eqref{eq:universal-mmm} and translation invariance,
\[
\langle t,z_0\rangle
=\langle\alpha_{\OP^2}^*\kappa_t,\xi\rangle,
\qquad t\in H^{24}(B\Spin(16);\Q).
\]
In particular, let $\Ahat_{24}\in H^{24}(B\Spin(16);\Q)$ denote the degree-$24$
component of the $\Ahat$-class. Since the ordinary family index of
$\xi$ vanishes by construction, we have
\[
\langle\Ahat_{24},z_0\rangle=0.
\]

We use the tangential structure $\theta$ defined below
to find a rational preimage $y\in\pi_8MT\theta\otimes\Q$ of
$z_0$ under a fixed-rank geometric transfer.
We then realize a positive multiple of the rational
Hurewicz image of $y$ by a smooth $\theta$-structured
family. The disk bundle caps associated to this $\theta$-family give a family $j$ of psc loops in the psc cobordism category. Comparing MMM pairings then produces the rational psc homotopy class $a_\Q$.

\subsection{Tangential structures and characteristic classes}
Define the tangential structure
\[
 \theta:B_8=B\Spin(8)\times B\Sp(1)\longrightarrow BO(8)
 \]
by forgetting the Spin lift and the quaternionic line. The associated Madsen--Tillmann spectrum is given by
\[
 MT\theta=MT\Spin(8)\wedge B\Sp(1)_+.
\]
For a smooth fiber bundle
\[
W \to X \xrightarrow{\pi} K
\]
with fiber $W$ a closed 8-manifold, a $\theta$-structure on $\pi$ consists of
\begin{itemize}
	\item A Spin structure on the vertical tangent bundle $T_\pi X$.
	\item A quaternionic line bundle $L\to X$ on the total space. 
\end{itemize}
These data determine a structure map $\ell_X:X\to B_8$
lifting the classifying map of $T_\pi X$ along $\theta$. We say such a fiber bundle is $\theta$-structured. 

Fix a closed 8-manifold $W$ with a chosen
$\theta$-structure $\ell:W\to B_8$ on the trivial
fibre bundle $W\to\mathrm{pt}$.
Let $\mathcal M^\theta(W,\ell)$ denote the moduli space
of smooth $\theta$-structured bundles whose fibre
structures are concordant to $(W,\ell)$.
A map $u:K\to\mathcal M^\theta(W,\ell)$ classifies
such a family up to concordance.
See \cite{GRWGuide} for details.

The fibrewise Pontryagin--Thom construction gives
\[
\alpha^{\mathrm{PT}}:\mathcal M^\theta(W,\ell)
\longrightarrow \Omega^\infty MT\theta.
\]
By translation in the infinite loop space, we have 
\[
\alpha:\mathcal M^\theta(W,\ell)
\longrightarrow \Omega^\infty_0MT\theta,
\qquad
\alpha(x)=\alpha^{\mathrm{PT}}(x)-\alpha^{\mathrm{PT}}(W, \ell).
\]
For $t\in H^{8+s}(B_8;\Q)$ with $s>0$,
let $\kappa_t\in H^s(\Omega^\infty MT\theta;\Q)$
be the universal MMM class defined as in
\eqref{eq:universal-mmm}, and write
\begin{equation}\label{eq:structured-mmm}
	\kappa_t^\theta
	=\alpha^*\kappa_t
	\in H^s(\mathcal M^\theta(W,\ell);\Q).
\end{equation}
For a $\theta$-structured family $\pi:X\to K$
classified by $u$, we have
\begin{equation}\label{eq:kappa-theta}
	u^*\kappa_t^\theta
	=\pi_!\ell_X^*t\in H^s(K;\Q).
\end{equation}
This is the $\theta$-structured version of
\cite[Section~1.1]{GRW14MMM}.

\subsection{The geometric transfer and its fixed-rank refinement}

In this subsection, we construct a fixed-rank transfer map $\mu_8$ whose rational image contains $z_0\otimes1$.
We first use the ordinary quaternionic projective geometric transfer in Spin bordism to calculate its stable
Pontryagin part.

Following the standard convention, we use right quaternionic lines
and right $\Hh$-linear maps. Let $\gamma^1\to B\Sp(1)$ be the universal quaternionic line,
and denote its pullback to $B_8$ by the same symbol.

\medskip

For a closed $m$-dimensional Spin manifold $B$ with a quaternionic
line bundle $L$, taking the associated projective bundle total space
defines a geometric transfer map, cf. \cite{KS93}
\[
\begin{aligned}
	\Psi :\Omega_m^\Spin(B\Sp(1))
	&\longrightarrow \Omega_{m+8}^\Spin,\\
	[B,L]&\mapsto
	[\mathbb P_{\Hh}(\underline{\Hh}^{2}\oplus L)].
\end{aligned}
\]
We show that
\[
\operatorname{im}(\Psi\otimes\Q)
=\ker\bigl(\Omega_*^\Spin\otimes\Q
\xrightarrow{\operatorname{ind}_{KO}\otimes\Q}
KO_*\otimes\Q\bigr).
\]
Indeed, for $n\ge2$, let $P_n$ be the total space of
$\mathbb P_{\Hh}(\underline{\Hh}^{2}\oplus\gamma^1)$
over $\HP^{n-2}$. 
By the rational structure of Spin bordism \cite{Stong68, KS93}, we have
\[
\Omega_*^\Spin\otimes\Q
=\Q\bigl[[K_4],[P_2],[P_3],\ldots\bigr],
\]
where $K_4$ is a Spin four-manifold of nonzero $\Ahat$-genus,
for example a K3 surface. Since $\operatorname{ind}_{KO}(K_4)\ne0$ and each $P_n$ has a psc metric, the polynomial presentation gives
\[
\ker(\operatorname{ind}_{KO}\otimes\Q)
=\bigl([P_n]\mid n\ge2\bigr).
\]
By construction, each $[P_n]$ lies in $\im\Psi$.
Since $\Psi$ is $\Omega_*^\Spin$-linear, its rationalized
image contains the ideal above. Conversely, all associated $\HP^2$-bundle
total spaces have psc metrics \cite[Section~3]{KS93}, so the rationalized
image is contained in the index kernel. This proves the claimed
image equality.

\medskip

Now we describe the dual of the transfer map in degree 24. Note that evaluation on stable tangent bundles and the quaternionic line gives perfect
pairings
\[
\begin{aligned}
	H^{24}(B\Spin;\Q)
	&\times\Omega_{24}^{\Spin}\otimes\Q\rightarrow\Q,\\
	H^{16}(B\Spin\times B\Sp(1);\Q)
	&\times\Omega_{16}^{\Spin}(B\Sp(1))\otimes\Q\rightarrow\Q.
\end{aligned}
\]
Let
$q:\mathbb P_{\Hh}(\underline{\Hh}^{2}\oplus\gamma^1)
\to B\Spin\times B\Sp(1)$, $\gamma_\infty$ the universal
stable Spin bundle, and $T_q$ the vertical tangent bundle of $q$.
Then the degree-$24$ dual transfer is
\[
\begin{aligned}
	\Psi_{24}^*:
	H^{24}(B\Spin;\Q)&\longrightarrow
	H^{16}(B\Spin\times B\Sp(1);\Q),\\
	t&\longmapsto
	q_!t\bigl(q^*\gamma_\infty
	\oplus T_q\bigr).
\end{aligned}
\]
Consequently, 
\[
\ker(\Psi_{24}^*)=\Q\Ahat_{24},
\]

\medskip

The above calculation controls the stable Pontryagin
classes. 
Since our purpose is retaining the characteristic data of the rank-$16$ Pontryagin--Thom class $z_0$, we also need to detect the Euler class.
The following construction gives a fixed-rank refinement
of the geometric transfer.
Consider the two fibre bundles over $B_8$:
\begin{equation}\label{eq:two-bundles-B8}
	\begin{aligned}
		q_{\mathcal E}:\mathcal E
		&=\mathbb P_{\Hh}(\underline{\Hh}^{2}\oplus\gamma^1)
		\longrightarrow B_8,\\
		q_{\mathcal D}:\mathcal D
		&=S\bigl(\Hom_{\Hh}(\gamma^1,\underline{\Hh}^{2})
		\oplus\R\bigr)
		\longrightarrow B_8.
	\end{aligned}
\end{equation}
Since $\Sp(1)$ is simply connected, the natural $\Sp(1)$-actions on fibres lift to the unique Spin frame bundles of $\HP^2$
and $S^8$.
Thus $\mathcal E$ and $\mathcal D$ carry the induced vertical Spin structures.

Let $\gamma_8\to B_8$ be the pullback of the universal rank 8 Spin
bundle from $B\Spin(8)$. Combining its pullback with the
vertical tangent bundles gives rank-16 Spin bundles
\[
\begin{aligned}
	q_{\mathcal E}^*\gamma_8\oplus T(\mathcal E/B_8)&\longrightarrow \mathcal E,\\
	q_{\mathcal D}^*\gamma_8\oplus T(\mathcal D/B_8)&\longrightarrow \mathcal D.
\end{aligned}
\]

Define
\begin{equation}\label{eq:dualgysin}
	\begin{aligned}
		\mu_8^*:H^{24}(B\Spin(16);\Q)&\longrightarrow H^{16}(B_8;\Q),\\
		t&\longmapsto(q_{\mathcal E})_!t(q_{\mathcal E}^*\gamma_8
		\oplus T(\mathcal E/B_8))\\
		&\quad-(q_{\mathcal D})_!t(q_{\mathcal D}^*\gamma_8
		\oplus T(\mathcal D/B_8)).
	\end{aligned}
\end{equation}
Under rational stable Hurewicz and the Thom isomorphisms,
its dual is denoted as
\begin{equation}\label{eq:transfer}
	\mu_8:\pi_8MT\theta\otimes\Q
	\longrightarrow\pi_8MT\Spin(16)\otimes\Q.
\end{equation}

\begin{proposition}\label{prop:transfer}
	We have
	\[
	\ker(\mu_8^*)=\Q\Ahat_{24},
	\]
	or equivalently,
	\[
	\im(\mu_8)=\ker\langle\Ahat_{24},-\rangle.
	\]
\end{proposition}

\begin{proof}

	The following classes form a basis of
	$H^{24}(B\Spin(16);\Q)$:
	\begin{equation}\label{eq:thirteen-coordinates}
		\begin{gathered}
			p_6,\ p_5p_1,\ p_4p_2,\ p_4p_1^2,\ p_3^2,\ p_3p_2p_1,\ p_3p_1^3,\\
			p_2^3,\ p_2^2p_1^2,\ p_2p_1^4,\ p_1^6,
			\qquad e_{16}p_2,\ e_{16}p_1^2.
		\end{gathered}
	\end{equation}
	The first eleven are stable Pontryagin polynomials; the remaining
	two involve the Euler class.
	For a stable Pontryagin polynomial $t\in H^{24}(B\Spin;\Q)$,
	the sphere term in \eqref{eq:dualgysin}
	\[
	(q_{\mathcal D})_!t\bigl(q_{\mathcal D}^*\gamma_8
	\oplus T(\mathcal D/B_8)\bigr)
	\]
	vanishes.
	Indeed, since
	$\mathcal D=S(\Hom_{\Hh}(\gamma^1,\underline{\Hh}^{2})\oplus\R)$,
	we have
	\[
	T(\mathcal D/B_8)\oplus\R
	\cong q_{\mathcal D}^*
	\bigl(\Hom_{\Hh}(\gamma^1,\underline{\Hh}^{2})\oplus\R\bigr).
	\]
	Hence
	\[
	\begin{aligned}
		&(q_{\mathcal D})_!
		t\bigl(q_{\mathcal D}^*\gamma_8
		\oplus T(\mathcal D/B_8)\bigr)\\
		&\qquad=
		t\bigl(\gamma_8\oplus
		\Hom_{\Hh}(\gamma^1,\underline{\Hh}^{2})\bigr)
		(q_{\mathcal D})_!(1)=0.
	\end{aligned}
	\]
	On the other hand, the projective term
	\[
	(q_{\mathcal E})_!
	t\bigl(q_{\mathcal E}^*\gamma_8
	\oplus T(\mathcal E/B_8)\bigr)
	\]
	is the restriction of $\Psi_{24}^*(t)$ to $B_8$.
	Since stabilization of $B\Spin(8)$ is injective
	on rational cohomology through degree $16$, we conclude that
	the kernel of $\mu_8^*$ on stable Pontryagin polynomials
	is $\Q\Ahat_{24}$.

	It remains to consider $e_{16}f$, where
	$f\in\Q\{p_2,p_1^2\}$.
	Let $\iota:B\Spin(8)\to B_8$ classify the trivial quaternionic line.
	After restriction, the two bundles have fibres $\HP^2$ and $S^8$.
	For either fibre $F$, we have
	$e(\gamma_8\boxplus TF)=e_8e(TF)$;
	therefore only the base-degree term of
	$f(\gamma_8\boxplus TF)$ contributes
	after multiplication by $e(TF)$. It follows that
	\begin{equation}\label{eq:eulerimages}
		\iota^*\mu_8^*(e_{16}f)
		=\bigl(\chi(\HP^2)-\chi(S^8)\bigr)e_8f
		=e_8f,\qquad f\in\Q\{p_2,p_1^2\}.
	\end{equation}
	For stable $t$, the class $\iota^*\mu_8^*(t)$ lies in
	$\Q[p_1,p_2,p_3,e_8^2]$.
	Since $H^*(B\Spin(8);\Q)=\Q[p_1,p_2,p_3,e_8]$, the two
	odd classes $e_8p_2$ and $e_8p_1^2$ are independent modulo the
	stable part. Thus no additional kernel arises from the Euler classes,
	and
	\[
	\ker(\mu_8^*)=\Q\Ahat_{24}.
	\]
	The image equality follows by rational duality.
\end{proof}

Recall the Pontryagin--Thom class $z_0$ satisfies
$\langle\Ahat_{24},z_0\rangle=0$.  Proposition~\ref{prop:transfer}
therefore allows us to choose
\begin{equation}\label{eq:preimage}
	y\in\pi_8MT\theta\otimes\Q,\qquad \mu_8(y)=z_0\otimes1.
\end{equation}

\subsection{The structured moduli space and a geometric representative}

To realize the rational Hurewicz image of $y$ geometrically in terms of actual fibre bundles, we use the
structured moduli spaces of Galatius--Randal-Williams.
For that purpose, we need a $\theta$-structured manifold $(W,\ell)$ such that the family Pontryagin--Thom
map on its moduli space is surjective on $H_8(-;\Q)$, as given by the following lemma. 

\begin{lemma}\label{lem:W-ell}
	There exist a closed $\theta$-structured 8-manifold
	$(W,\ell)$, such that the structure map $\ell:W\to B_8$ is $4$-connected and the stable $\theta$-genus satisfies $\bar g^\theta(W,\ell)\ge28$.
\end{lemma}
\begin{proof}
	Let
	\[
	W=\HP^2\#\HP^2\#^{28}(S^4\times S^4).
	\]
	Equip $W$ with its Spin structure, and a quaternionic line
	bundle $L$ obtained by gluing the tautological line bundle
	on the second copy of $\HP^2$ to the trivial line bundles
	on the remaining summands.
	Let $\ell:W\to B_8$ be the resulting structure map.
	
	Since both $W$ and $B_8$ are $3$-connected, it suffices
	to show that
	\[
	\ell_*:H_4(W;\Z)\longrightarrow H_4(B_8;\Z)
	\]
	is surjective. Recall that $H^4(B_8;\Z)$ is freely generated
	by $p_1/2$ and $c_2$.
	Since $p_1(T\HP^2)/2$ generates $H^4(\HP^2;\Z)$,
	and $c_2(L)$ vanishes on the first $\HP^2$ summand
	and generates the degree-$4$ cohomology of the second,
	the two summand generators map to a basis of
	$H_4(B_8;\Z)$. Thus $\ell$ is $4$-connected
	by the Hurewicz theorem.

	The $28$ copies of $S^4\times S^4$ in
	the connected sum contribute $28$ disjoint embedded copies of
	$(S^4\times S^4)\setminus\operatorname{int}(D^8)$, 
	each with the unique Spin structure
	and trivial quaternionic line bundle.
	Hence the
	$\theta$-structures on the copies are admissible, see \cite[Remark~2.11]{GRW17}.
	It follows that
	\[
	\bar g^\theta(W,\ell)\ge g^\theta(W,\ell)\ge28.
	\]
	This gives the desired $\theta$-structured manifold. 
\end{proof}

\medskip

Fix $(W,\ell)$ as in Lemma~\ref{lem:W-ell}.
Consider the structured moduli space $\mathcal M_4^\theta(W,\ell)$,
which classifies smooth $\theta$-structured bundles with fibre
structures concordant to $(W,\ell)$. The subscript $4$ records
that the fibre structure maps are $4$-connected, as follows
from their concordance to $\ell$.
Recall the normalized family Pontryagin--Thom map
defined above
\begin{equation}\label{eq:grw}
	\alpha:\mathcal M_4^\theta(W,\ell)
	\longrightarrow\Omega^\infty_0MT\theta. 
\end{equation}
Since $B_8=B\Spin(8)\times B\Sp(1)$ is simply connected, Galatius and Randal-Williams' results imply that $\alpha$
induces an isomorphism on $H_i(-;\Q)$,  whenever the stable $\theta$-genus $\bar g^\theta(W,\ell)$ is no less than $3i+4$. See \cite{GRW17} for details.
In our case of interest, $\bar g^\theta(W,\ell)\ge28$ suffices for the surjectivity
\[
	\alpha_*:
	H_8(\mathcal M_4^\theta(W,\ell); \Q)
	\twoheadrightarrow
	H_8(\Omega^\infty_0MT\theta; \Q). 
\]
Combining with the fact that the Thom orientation $MSO\to H\Z$ induces a surjection on rational homology theories, we arrive at
\begin{equation}\label{eq:surjections}
	\Omega_8^{SO}(\mathcal M_4^\theta(W,\ell))\otimes\Q
	\twoheadrightarrow
	H_8(\mathcal M_4^\theta(W,\ell);\Q)
	\overset{\alpha_*}{\twoheadrightarrow}
	H_8(\Omega^\infty_0MT\theta; \Q)
	\xleftarrow{\operatorname{hur}_{\Q}}
	\pi_8MT\theta\otimes\Q. 
\end{equation}
After clearing denominators, we may choose a positive integer $r$, a possibly disconnected closed oriented 8-manifold $K$, and a map
$u:K\to\mathcal M_4^\theta(W,\ell)$ such that
\begin{equation}\label{eq:cycle}
	\frac1r[K,u]\longmapsto\frac1r u_*[K]
	\overset{\alpha_*}{\longmapsto}
	\operatorname{hur}_{\Q}(y)
	\overset{\operatorname{hur}_{\Q}}{\longmapsfrom}y.
\end{equation}
Concretely, $u$ classifies a smooth $\theta$-structured bundle
\begin{equation}\label{eq:actual-family}
	\pi:X\longrightarrow K
\end{equation}
with structure map $\ell_X:X\to B_8$ satisfying
$\ell_X^*\gamma_8\cong T(X/K)$.
For each $k\in K$, the restricted structured fibre
$(X_k,\ell_X|_{X_k})$ is concordant to $(W,\ell)$.

\medskip
For $t\in H^{16}(B_8;\Q)$, combining
\eqref{eq:structured-mmm}, \eqref{eq:kappa-theta},
and \eqref{eq:cycle} gives
\begin{equation}\label{eq:realized-pairing}
	\begin{aligned}
		\langle\ell_X^*t,[X]\rangle
		&=\langle\kappa_t^\theta,u_*[K]\rangle\\
		&=\langle\kappa_t,\alpha_*u_*[K]\rangle\\
		&=r\langle\kappa_t,\operatorname{hur}_{\Q}(y)\rangle\\
		&=r\langle t,y\rangle.
	\end{aligned}
\end{equation}
The orientation of $X$ is induced by those of
$T(X/K)$ and $K$. The last equality follows from
the definition of $\kappa_t$.

\subsection{Cobordism categories and a rational psc lift}\label{subsec:spherical}

We now use the $\theta$-structured family $X\to K$
constructed above to produce a family of psc loops.
We will identify its forgetful image with the difference of two closed gluings, and then
use the resulting characteristic pairings to recover
$z_0\otimes1$ in rational homotopy.

\medskip

We use the cobordism categories of Ebert--Randal-Williams
\cite[Definitions~3.1.1 and 3.2.1--3.2.2]{ERW22}.
The ordinary Spin cobordism category $\cC_{\Spin(16)}$
has closed 15-dimensional Spin manifolds
as objects, and collared 16-dimensional Spin cobordisms as morphisms.
Its subcategory $\cC_{\Spin(16)}^2$ consists of cobordisms
for which the outgoing boundary inclusion 
is $2$-connected.
The psc cobordism category $\cP_{\Spin(16)}^{2,\rst}$
equips objects with psc metrics and morphisms with right stable
psc metrics which are products near collars.
Here a morphism metric is right stable if gluing it to any
following cobordism induces a weak homotopy equivalence
on the relevant spaces of psc metrics \cite[Definition~2.3.1]{ERW22}.

These categories give maps of classifying spaces
\begin{equation}\label{eq:forgetful-categories}
	f: B\cP_{\Spin(16)}^{2,\rst}
	\xrightarrow{f_2} B\cC_{\Spin(16)}^2
	\xrightarrow{i} B\cC_{\Spin(16)}. 
\end{equation}
Here $f_2$ forgets the psc metrics, while the inclusion $i$ induces a weak equivalence by
\cite[Theorem~3.1.3]{ERW22}.

\medskip

We now decompose the two fibre bundles in \eqref{eq:two-bundles-B8} as equivariant gluings of disk bundle caps.
Recall
the $\Sp(1)$-action on $\HP^2=\mathbb P_{\Hh}(\Hh^2\oplus\Hh)$ is given by
\[
a\cdot[v,w]=[v,aw],
\qquad a\in\Sp(1), 
\quad
v\in\Hh^2,\quad w\in\Hh.
\]
Let $\gamma^1\to\HP^1$ be the tautological quaternionic line.
The disk bundle
\[
D\bigl(\Hom_{\Hh}(\gamma^1,\underline{\Hh})\bigr)
\longrightarrow\HP^1
\]
is identified with the region $\{\,[v,w]:\|w\|\le\|v\|\,\} \subset \mathbb P_{\Hh}(\Hh^2\oplus\Hh)$ by the graph map
\[
(\lambda,\phi)\longmapsto[v,\phi(v)],
\qquad 0\ne v\in\lambda, 
\quad \lambda \in \HP^1. 
\]
Attaching the cap $D(\Hh^2)=\{\,[x,1]:\|x\|\le 1\,\}$
via the boundary identification
\[
g:
S\bigl(\Hom_{\Hh}(\gamma^1,\underline{\Hh})\bigr)
\longrightarrow S(\Hh^2),
\qquad
g(\lambda,\phi)=\phi^{-1}(1),
\]
we arrive at an $\Sp(1)$-equivariant splitting 
\begin{equation}\label{eq:split}
	\HP^2=
	D\bigl(\Hom_{\Hh}(\gamma^1,\underline{\Hh})\bigr)
	\cup_g D(\Hh^2)^{\op}.
\end{equation}
Indeed, the actions on the two pieces are given by
$(\lambda,\phi)\mapsto(\lambda,a\phi)$ and
$[x, 1]\mapsto [xa^{-1},1 ]$, respectively, and the gluing map is equivariant since
$g(\lambda,a\phi)=g(\lambda,\phi)a^{-1}$.
See \cite[(3.20)]{EW24} for the general case. 
Here we view the two pieces as outgoing cobordisms from the empty set to their common boundary, and the superscript op denotes cobordism reversal.

\medskip

For the $\theta$-structured family $\pi:X\to K$
in \eqref{eq:actual-family}, let $P\to X$ be the
principal $\Sp(1)$-bundle over $X$, with its associated quaternionic line denoted $L$.
Form the disk bundles
\begin{equation}\label{eq:two-caps}
	\begin{aligned}
		Q_0&=P\times_{\Sp(1)}
		D\bigl(\Hom_{\Hh}(\gamma^1,\underline{\Hh})\bigr)
		\longrightarrow \HP^1\times X\longrightarrow X,\\
		Q_1&=P\times_{\Sp(1)}D(\Hh^2)
		\cong D\bigl(\Hom_{\Hh}(L,\underline{\Hh}^{\,2})\bigr)
		\longrightarrow X.
	\end{aligned}
\end{equation}
The map induced by $g$ identifies their boundaries with
$M=P\times_{\Sp(1)}S^7$.
Gluing along $M$ gives two closed bundles over $X$:
\begin{equation}\label{eq:two-gluings}
	\begin{aligned}
		E&=Q_0\cup_M Q_1^{\op}
		=\mathbb P_{\Hh}(\underline{\Hh}^{\,2}\oplus L)
		\xrightarrow{q_E}X,\\
		D&=Q_1\cup_M Q_1^{\op}
		=S\bigl(\Hom_{\Hh}(L,\underline{\Hh}^{\,2})\oplus\R\bigr)
		\xrightarrow{q_D}X,
	\end{aligned}
\end{equation}
with fibers $\HP^2$ and $S^8$, respectively.
Since $\ell_X^*\gamma_8\cong T(X/K)$ and
$\ell_X^*\gamma^1\cong L$, the bundles $E$ and $D$
in \eqref{eq:two-gluings} are the pullbacks of
$\mathcal E$ and $\mathcal D$ along $\ell_X$, with their
induced vertical Spin structures.
The composites $\pi\circ q_E$
and $\pi\circ q_D$ are families of closed 16-manifolds over
$K$, with Spin splittings
\[
\begin{aligned}
	T(E/K)&\cong q_E^*T(X/K)\oplus T(E/X),\\
	T(D/K)&\cong q_D^*T(X/K)\oplus T(D/X).
\end{aligned}
\]
For $t\in H^{24}(B\Spin(16);\Q)$, 
\eqref{eq:dualgysin} gives
\[
\ell_X^*(\mu_8^*t)
=(q_E)_!t(T(E/K))-(q_D)_!t(T(D/K)).
\]
Applying $\pi_!$ gives the difference of the two MMM classes.

\medskip

For each $k\in K$, the disk bundle caps in \eqref{eq:two-caps}
give Spin cobordisms
$Q_{0,k},Q_{1,k}:\varnothing\to M_k$ with $2$-connected outgoing boundary inclusions, 
since the fibres are the disks $D^4$ and $D^8$.
To construct a loop family in the psc cobordism, we seek smooth
families of right stable psc metrics $m_0,m_1$
on these caps which are products near boundary collars,
such that
\begin{equation}\label{eq:two-morphisms}
	(Q_{i,k},m_i|_{Q_{i,k}}):
	\varnothing\longrightarrow (M_k,m_i|_{M_k}),
	\qquad i=0,1,
\end{equation}
are morphisms in $\cP_{\Spin(16)}^{2,\rst}$.

\begin{lemma}\label{lem:caps}
	Such metrics $m_0,m_1$ exist with $m_0|_M=m_1|_M$.
\end{lemma}

\begin{proof}
	Let $b_X$ be a metric 
	on $T(X/K)$, and write $g_{S^n}$ for the unit round metric.	
	Choose a connection $\omega$ on $P\to X$, inducing
	connections on the associated bundles $M,Q_0,Q_1$ over $X$.
	By \cite[(3.13)]{EW24} and compactness of $X$, we may choose a sufficiently
	small $\delta>0$ such that the connection metric
	\[
	a=b_X\oplus_\omega\delta^2g_{S^7}
	\]
	on $T(M/K)$ is fibrewise psc.
	For $n=4,8$, choose a torpedo metric $g_{\mathrm{tor}}^n$
	of radius $\delta$ on $D^n$ as in
	\cite[Section~2.1]{EW24}.
	On $T(Q_1/K)$, set
	\[
	m_1=b_X\oplus_\omega g_{\mathrm{tor}}^8.
	\]
	For $Q_0$, let $\omega_{\mathrm H}$ be the metric connection
	on $\Hom_{\Hh}(\gamma^1,\underline{\Hh})\to\HP^1$
	induced by the standard Hopf connection on
	$S^7\to\HP^1$.
	The splitting induced by $\omega_{\mathrm H}$ is
	$\Sp(1)$-equivariant and hence descends to a splitting of
	\[
	0\longrightarrow T(Q_0/(\HP^1\times X))
	\longrightarrow T(Q_0/X)
	\longrightarrow \operatorname{pr}_{\HP^1}^*T\HP^1
	\longrightarrow0.
	\]
	Let $g_{\HP^1}$ be the quotient metric of $g_{S^7}$
	under the Hopf fibration.
	The corresponding connection metric on $T(Q_0/X)$ is
	\[
	\delta^2g_{\HP^1}
	\oplus_{\omega_{\mathrm H}}g_{\mathrm{tor}}^4.
	\]
	Using the splitting induced by $\omega$, we then define
	\[
	m_0=b_X\oplus_\omega
	\bigl(\delta^2g_{\HP^1}
	\oplus_{\omega_{\mathrm H}}g_{\mathrm{tor}}^4\bigr)
	\]
	on $T(Q_0/K)$.
	Under the boundary identification $g$, we have
	\[
	\delta^2g_{\HP^1}
	\oplus_{\omega_{\mathrm H}}\delta^2g_{S^3}
	=g^*(\delta^2g_{S^7}).
	\]
	Since $g$ is $\Sp(1)$-equivariant, it
	also preserves the splittings induced
	by $\omega$ on the associated boundary bundles.
	Hence
	\[
	m_0|_M=m_1|_M=a.
	\]
	The torpedo metrics give common product collars
	with metric $dt^2+a$.

	For each $k\in K$, the metrics $m_0|_{Q_{0,k}}$ and
	$m_1|_{Q_{1,k}}$ are connection torpedo metrics
	on disk bundles of ranks $4$ and $8$ over
	\[
	\bigl(\HP^1\times X_k,
	\delta^2g_{\HP^1}\oplus b_X|_{X_k}\bigr)
	\quad\text{and}\quad
	\bigl(X_k,b_X|_{X_k}\bigr),
	\]
	respectively, by construction and \cite[Lemma~3.7]{EW24}.
	The O'Neill formula
	and the properties of torpedo metrics
	\cite[Section~2.1]{EW24} give
	\[
	\scal(m_i|_{Q_{i,k}})
	\ge \min_{M_k}\scal(a|_{M_k})>0,
	\qquad i=0,1.
	\]
	Since the disk bundles have ranks $4$ and $8$,
	their connection torpedo metrics are right stable
	by \cite[Theorem~3.15]{EW24}.

	In the decomposition \eqref{eq:split},
	$g$ preserves boundary orientations and admits an
	equivariant Spin lift.
	Together with the fixed Spin structure on $T(X/K)$,
	this identifies the boundary Spin structures of $Q_0$
	and $Q_1$.
	Choosing common boundary Spin representatives and embeddings,
	the two morphisms in \eqref{eq:two-morphisms} have the same target.
\end{proof}

\medskip
In the geometric realization $B\cP_{\Spin(16)}^{2,\rst}$,
the two morphisms in \eqref{eq:two-morphisms}
determine paths with a common endpoint:
\begin{equation}\label{eq:two-paths}
	q_i(k):\varnothing\rightsquigarrow(M_k,a|_{M_k}),
	\qquad i=0,1.
\end{equation}
The concatenation gives a based loop
\[
j(k):=q_0(k)*\overline{q_1(k)}
:\varnothing\rightsquigarrow\varnothing.
\]
Varying these loops with $k\in K$ gives a family
\[
K\xrightarrow{j}\Omega_\varnothing
B\cP_{\Spin(16)}^{2,\rst}
\xrightarrow{\Omega f}
\Omega_\varnothing B\cC_{\Spin(16)}
\xrightarrow[\simeq]{\operatorname{PT}}
\Omega^\infty MT\Spin(16),
\]
where the last map is the Pontryagin--Thom equivalence
of \cite{GMTW09}. 

\medskip

Next we relate the loop family $j$ to the two closed gluings $E$ and $D$. Write
\[
\begin{aligned}
	\alpha_E:
	K\longrightarrow\Omega_\varnothing\Bord, \qquad
	\alpha_D:
	K\longrightarrow\Omega_\varnothing\Bord,
\end{aligned}
\]
for the family Pontryagin--Thom maps of the smooth families~\eqref{eq:two-gluings} 
\[
\begin{aligned}
	\pi_E: E\xrightarrow{q_E} X
	\xrightarrow{\pi} K,\qquad
	\pi_D: D\xrightarrow{q_D} X
	\xrightarrow{\pi} K,
\end{aligned}
\]
of closed Spin 16-manifolds over $K$. 

\begin{lemma}\label{lem:compare}
	We have the following homotopy
	\begin{equation}\label{eq:loopcompare}
		(\Omega f)\circ j\simeq\alpha_E*\overline{\alpha_D}
	\end{equation}
	as maps $K\rightarrow\Omega_\varnothing\Bord$.
\end{lemma}

\begin{proof}
	By \eqref{eq:two-paths}, the paths $f\circ q_0(k)$
	and $f\circ q_1(k)$ have the same endpoint $M_k$.
	Thus
	\[
	(\Omega f)(j(k))\simeq
	(f\circ q_0(k))*\overline{(f\circ q_1(k))}.
	\]
	Meanwhile, the cobordism
	$Q_{1,k}^{\op}:M_k\rightarrow\varnothing$ gives rise to a path
	\[
	(f\circ q_1(k))^{\op}:M_k\rightsquigarrow\varnothing
	\]
	in $\Bord$. 
	By \eqref{eq:two-gluings}, composing $Q_{0,k}$ and $Q_{1,k}$
	with $Q_{1,k}^{\op}$ gives the closed Spin manifolds
	$E_k$ and $D_k$, respectively.
	The Pontryagin--Thom construction then gives
	\[
	\begin{aligned}
		\alpha_E(k)&\simeq(f\circ q_0(k))*(f\circ q_1(k))^{\op},\\
		\alpha_D(k)&\simeq(f\circ q_1(k))*(f\circ q_1(k))^{\op}.
	\end{aligned}
	\]
	The composition of $2$-simplices in the nerve, together with
	the Pontryagin--Thom description of closed morphisms, give
	\[
	\begin{aligned}
		\alpha_E(k)*\overline{\alpha_D(k)}
		&\simeq(f\circ q_0(k))*(f\circ q_1(k))^{\op}
		*\overline{(f\circ q_1(k))^{\op}}
		*\overline{(f\circ q_1(k))}\\
		&\simeq(f\circ q_0(k))*\overline{(f\circ q_1(k))}
		\simeq(\Omega f)(j(k)).
	\end{aligned}
	\]
	These homotopies vary continuously over $K$, with the common
	collars and Spin identifications chosen in the proof of Lemma~\ref{lem:caps}.
\end{proof}

For $t\in H^{16+s}(B\Spin(16);\Q)$ with $s>0$,
we use the Pontryagin--Thom equivalence to identify
$\kappa_t$ with its pullback to
$H^s(\Omega_\varnothing\Bord;\Q)$.
Applying \eqref{eq:family-mmm} to the two closed families gives
\begin{equation}\label{eq:mmm}
	\alpha_E^*\kappa_t
	=(\pi_E)_!t(T(E/K)),
	\qquad
	\alpha_D^*\kappa_t
	=(\pi_D)_!t(T(D/K)),
\end{equation}
in $H^s(K;\Q)$.
Since $\operatorname{PT}$ respects loop multiplication,
$\kappa_t$ is primitive for loop concatenation.
For $t\in H^{24}(B\Spin(16);\Q)$, we obtain
\begin{equation}\label{eq:classcompare}
	\begin{aligned}
		j^*(\Omega f)^*\kappa_t
		&=\alpha_E^*\kappa_t-\alpha_D^*\kappa_t\\
		&=(\pi_E)_!t(T(E/K))-(\pi_D)_!t(T(D/K))\\
		&=\pi_!\ell_X^*(\mu_8^*t)
		=u^*\kappa^\theta_{\mu_8^*t},
	\end{aligned}
\end{equation}
which follows from \eqref{eq:loopcompare}, \eqref{eq:dualgysin}, and 
 \eqref{eq:kappa-theta}.

On each connected component of $K$, we may translate $j$
by the inverse of a fixed loop in its image, so that it takes
values in the component of the constant loop. Since this
does not change positive degree primitive pairings, we continue
to write
\[
j:K\longrightarrow\Omega_0\Bpsc
\]
for the resulting map. This choice of convention enables us to apply the result of Milnor--Moore later.
For $t\in H^{24}(B\Spin(16);\Q)$, we pair
$(\Omega f)^*\kappa_t$ with $j_*[K]$:
\begin{equation}\label{eq:pairing-chain}
	\begin{aligned}
		\langle(\Omega f)^*\kappa_t,j_*[K]\rangle
		=\langle\kappa^\theta_{\mu_8^*t},u_*[K]\rangle
		=r\langle\mu_8^*t,y\rangle
		=r\langle t,z_0\rangle.
	\end{aligned}
\end{equation}
Here we use \eqref{eq:classcompare},
\eqref{eq:realized-pairing}, and \eqref{eq:preimage}.
The chain holds for every $t$,
including the Euler coordinates $e_{16}p_2$ and
$e_{16}p_1^2$ in \eqref{eq:thirteen-coordinates}.

\medskip

As a homology class, $j_*[K]$ is not necessarily spherical. We now construct a rational
homotopy class
\[
a_\Q\in\pi_8(\Omega_0\Bpsc)\otimes\Q
\]
whose Hurewicz image has the same primitive pairings as
$r^{-1}j_*[K]$. We use the following result of
Milnor--Moore.

\begin{lemma}\label{lem:primitive}
	Let $Z$ be a based space of CW type, and let
	$I$ be the positive-degree ideal of $H_*(\Omega_0Z;\Q)$ under the
	Pontryagin product. For $s>0$, the composite
	\[
	\pi_s(\Omega_0Z)\otimes\Q
	\xrightarrow{\operatorname{hur}_{\Q}}
	H_s(\Omega_0Z;\Q)
	\longrightarrow H_s(\Omega_0Z;\Q)/(I^2)_s
	\]
	is surjective, where $(I^2)_s$ is the degree-$s$ part of $I^2$.
\end{lemma}

\begin{proof}
	The connected loop space $\Omega_0Z$ is a homotopy associative
	$H$-space. By \cite[Appendix]{MM65}, its rational
	Pontryagin homology algebra is generated by rational Hurewicz images.
	Modulo $I^2$, only products of length one remain, proving the
	surjectivity.
\end{proof}

Apply Lemma~\ref{lem:primitive} to $Z=\Bpsc$ and $s=8$. We obtain
\begin{equation}\label{eq:spherical-class}
	a_\Q \in\pi_8(\Omega_0\Bpsc)\otimes\Q,
	\qquad
	r^{-1} j_*[K]-\operatorname{hur}_{\Q}(a_\Q)\in I^2.
\end{equation}
For $t\in H^{24}(B\Spin(16);\Q)$, the class
$(\Omega f)^*\kappa_t$ is primitive because
$\Omega f$ preserves Pontryagin products.
It therefore pairs trivially with $I^2$.
By the naturality of $\operatorname{hur}_{\Q}$, we have
\[
\begin{aligned}
	\langle\kappa_t,
	\operatorname{hur}_{\Q}((\Omega f)_*a_\Q)\rangle
	&=\langle\kappa_t,
	(\Omega f)_*\operatorname{hur}_{\Q}(a_\Q)\rangle\\
	&=\langle\kappa_t,
	r^{-1}(\Omega f)_*j_*[K]\rangle\\
	&=\langle t,z_0\rangle.
\end{aligned}
\] 
Rational stable Hurewicz and the Thom isomorphism show
that these pairings determine the class in
$\pi_8MT\Spin(16)\otimes\Q$.
Hence
\begin{equation}\label{eq:spherical-image}
	\operatorname{PT}_*(\Omega f)_*a_\Q=z_0\otimes1.
\end{equation}
The following diagram records our construction above.

\begin{scriptsize}
	\[\begin{tikzcd}[column sep=2.5em]
		& {\pi_8\bigl(\Omega_0 B\mathcal{P}_{\mathrm{Spin}(16)}^{2,\mathrm{rst}}\bigr)\otimes\mathbb{Q}} & {\pi_8\bigl(\Omega_0 B\mathcal{C}_{\mathrm{Spin}(16)}\bigr)\otimes\mathbb{Q}} & {\pi_8 MT\mathrm{Spin}(16)\otimes\mathbb{Q}} & \\
		{H_8(K;\mathbb{Q})} & {H_8\bigl(\Omega_0 B\mathcal{P}_{\mathrm{Spin}(16)}^{2,\mathrm{rst}};\mathbb{Q}\bigr)} & {H_8\bigl(\Omega_0B\mathcal{C}_{\mathrm{Spin}(16)};\mathbb{Q}\bigr)} & {H_{8}(MT\mathrm{Spin}(16);\mathbb{Q})} & {H_{24}(B\mathrm{Spin}(16);\mathbb{Q})} \\
		&& {H^8\bigl(\Omega_0B\mathcal{C}_{\mathrm{Spin}(16)};\mathbb{Q}\bigr)} & {H^{8}(MT\mathrm{Spin}(16);\mathbb{Q})} & {H^{24}(B\mathrm{Spin}(16);\mathbb{Q})}
		\arrow["{(\Omega f)_*}", from=1-2, to=1-3]
		\arrow["{\operatorname{hur}_{\mathbb Q}}"', from=1-2, to=2-2]
		\arrow["{\mathrm{PT}_*}", hook', two heads, from=1-3, to=1-4]
		\arrow["{\operatorname{hur}_{\mathbb Q}}"', from=1-3, to=2-3]
		\arrow["{\operatorname{hur}_{\mathbb Q}}"', hook', two heads, from=1-4, to=2-4]
		\arrow["{j_*}", from=2-1, to=2-2]
		\arrow["{(\Omega f)_*}", from=2-2, to=2-3]
		\arrow["{\sigma_*\circ\mathrm{PT}_*}", from=2-3, to=2-4]
		\arrow["{\langle-,-\rangle}"{description}, draw=none, from=2-3, to=3-3]
		\arrow["{\mathrm{Thom}}", hook', two heads, from=2-4, to=2-5]
		\arrow["{\langle-,-\rangle}"{description}, draw=none, from=2-4, to=3-4]
		\arrow["{\langle-,-\rangle}"{description}, draw=none, from=2-5, to=3-5]
		\arrow["{\mathrm{PT}^*\circ\sigma}"', from=3-4, to=3-3]
		\arrow[hook', two heads, from=3-5, to=3-4]
	\end{tikzcd}\]
	\end{scriptsize}

	Thus $a_\Q$ is a rational psc lift of the
	Pontryagin--Thom class $z_0$ of the original family.

\section{Fibrewise positive scalar curvature on the Cayley family}\label{sec:based}

Starting from the punctured bundle, we now construct smooth fibrewise psc metrics on a
positive degree pullback of the original Cayley bundle represented
by $\xi$. 

\subsection{The punctured bundle}
Choose a parametrized embedded disk
$D^{16}\subset\OP^2$ with a boundary collar, and set
\[
N=\OP^2\setminus\operatorname{int}D^{16}\simeq S^8,\qquad
\partial N=S^{15}.
\]
The unique Spin structure on $N$ is compatible with
the standard Spin structure on the disk.

Consider a Spin preserving diffeomorphism $f$ of $N$
satisfying:
\begin{itemize}
	\item $f$ is the identity on the boundary collar.
	\item Extending by identity over $D^{16}$, $f$ extends to a diffeomorphism in
	$\widehat{\Diff}_0(\OP^2)$.
\end{itemize}
Denote the group of such diffeomorphism by $\widehat{\Diff}_\partial(N)$, with the natural capping map
\[
\operatorname{cap}:\widehat{\Diff}_\partial(N)
\longrightarrow\widehat{\Diff}_0(\OP^2).
\]

\begin{lemma}\label{lem:disk}
	The capping map induces a surjection
	\[
	(B\operatorname{cap})_*: 
	\pi_8B\widehat{\Diff}_\partial(N)\longrightarrow\pi_8\Bdiff
	\]
	given by gluing in the trivial disk
	bundle $S^8\times D^{16}$.
	In particular, $\xi$ has an integral lift $\xi_\partial$. 
\end{lemma}
\begin{proof}
	The disk restriction fibration has base homotopy equivalent to the oriented frame bundle of $T\OP^2$. Its Spin cover has the same higher homotopy groups. The frame sequence contains
	\[
	\pi_8(\OP^2)\longrightarrow\pi_7\Spin(16)
	\longrightarrow\pi_7\Fr^\Spin(T\OP^2)
	\longrightarrow\pi_7(\OP^2).
	\]
	The first two groups are $\Z$ and the last is zero. Complexification takes a real Bott generator on $S^8$ to a complex generator; since $\ph_2=-p_2/6$ there, its $p_2$ is $\pm6$ times the integral generator. The equality $p_2(T\OP^2)=6u$ makes the first arrow multiplication by $\pm1$. Thus $\pi_7\Fr^\Spin(T\OP^2)=0$, and the disk-restriction exact sequence proves the assertion.
\end{proof}

Consider the family Pontryagin--Thom
map of the capped $N$-bundles:
\[
\alpha_N:=\alpha_{\OP^2}\circ B\operatorname{cap}:
B\widehat{\Diff}_\partial(N)
\longrightarrow\Omega^\infty MT\Spin(16),
\]
as in \cite[proof of Theorem~F]{ERW22}.
Fix a lift $\xi_\partial$ of $\xi$ given by Lemma~\ref{lem:disk}.
Since $(B\operatorname{cap})_*\xi_\partial=\xi$ and
$\alpha_N(N)=\alpha_{\OP^2}(\OP^2)$, we may equivalently express the normalized class \eqref{eq:z} in terms of $\xi_\partial$ as follows
\begin{equation}\label{eq:z2}
	z_0=
	\bigl(\alpha_N(-)-\alpha_N(N)\bigr)_*\xi_\partial.
\end{equation}

The remaining metric lifting problem is to find a positive
integer $k$ such that the degree-$k$ pullback of the smooth $N$-bundle $\xi_\partial$ admits fibrewise psc metrics with a fixed
product collar.
Gluing in a fixed torpedo metric on
$D^{16}$ with the same boundary value, 
we obtain fibrewise psc metrics on the closed Cayley family
represented by $k\xi$.

\subsection{The psc metric lifting problem}\label{subsec:based-model}

In this subsection, we discuss the psc metric lifting obstruction on the punctured bundle, and compare the family Pontryagin--Thom class
of $\xi_\partial$ with the rational psc homotopy class $a_\Q$ constructed in Section~\ref{sec:metrics}. 

\medskip

Let $\varphi:\R^{16}\hookrightarrow\OP^2$ extend the chosen
disk embedding $D^{16}\hookrightarrow\OP^2$.
Choose a torpedo metric $g_{\mathrm{tor}}$ on $\R^{16}$
and write its product form near $\partial D^{16}=S^{15}$
as $g_{\mathrm{tor}}=dt^2+g_\partial$.
By \cite{ERW22, Che04},
the inclusion
\[
\cR^+(\OP^2,\varphi)
:=
\left\{
g\in\cR^+(\OP^2):
\varphi^*g=g_{\mathrm{tor}}
\text{ near }D^{16}
\right\}
\xhookrightarrow{}
\cR^+(\OP^2)
\]
is a weak homotopy equivalence. 
Then by surjectivity on path components, we may choose
$g\in\cR^+(\OP^2,\varphi)$ in the path component of the standard metric 
$g_{\mathrm{std}}$ on $\OP^2$. We fix a path from $g$ to
$g_{\mathrm{std}}$ in $\cR^+(\OP^2)$.

Following \cite[Section~2.1]{ERW22}, let
\[
\cR^+(N)_{g_\partial}
=
\left\{
m\in\cR^+(N):
m=dt^2+g_\partial
\text{ near }\partial N
\right\}.
\]
By construction, there is a restriction map
\[
\begin{aligned}
	\cR^+(\OP^2,\varphi)
	\longrightarrow\cR^+(N)_{g_\partial}.
\end{aligned}
\]
Apply the restriction to $g$, we have
\[
g=h\cup g_{\mathrm{tor}},
\qquad h:=g|_N\in\cR^+(N)_{g_\partial}.
\]

\medskip

The group $\widehat{\Diff}_\partial(N)$ acts on
$\cR^+(N)_{g_\partial}$ through its underlying
diffeomorphisms. Explicitly, a Spin lift
$\widehat\varphi$ acts by
$(\varphi^{-1})^*m$, where $\varphi$ is the underlying
diffeomorphism, which is the identity on the boundary collar. 
This gives the associated bundle
\begin{equation}\label{eq:metric-borel}
	\cR^+(N)_{g_\partial}
	\longrightarrow
	E\widehat{\Diff}_\partial(N)
	\times_{\widehat{\Diff}_\partial(N)}
	\cR^+(N)_{g_\partial}
	\longrightarrow
	B\widehat{\Diff}_\partial(N).
\end{equation}
Pulling back along a classifying map
$S^8\to B\widehat{\Diff}_\partial(N)$ representing
$\xi_\partial$ gives the bundle of psc metric spaces
of the specified $N$-family. A section of this bundle
is a continuous family of psc metrics with the fixed
boundary condition.
Choose a basepoint $e\in E\widehat{\Diff}_\partial(N)$ over
the basepoint $N\in B\widehat{\Diff}_\partial(N)$.
Then $[e,h]$ is a basepoint of the associated bundle.

The connecting homomorphism is
\[
\delta_h:
\pi_8(B\widehat{\Diff}_\partial(N),*)
\longrightarrow
\pi_7(\cR^+(N)_{g_\partial},h).
\]
For a based clutching map
$S^7\to\widehat{\Diff}_\partial(N)$ representing
$\xi_\partial$, the class $\delta_h(\xi_\partial)$
is represented by its composition with the orbit map
$\widehat\varphi\mapsto\widehat\varphi\cdot h$.
Since the base is $S^8$, this obstruction vanishes
precisely when the metric bundle admits a continuous
section taking the value $h$ at the basepoint. 
It therefore suffices to prove that
$\delta_h(\xi_\partial)$ is torsion.

\medskip

The underlying diffeomorphisms give a map
\[
B\widehat{\Diff}_\partial(N)
\longrightarrow B\Diff_\partial(N).
\]
The morphism space
$\cC_{\Spin(16)}^2(\varnothing,S^{15})_N$
consists of cobordisms from $\varnothing$ to $S^{15}$
whose underlying manifolds are diffeomorphic to $N$
relative to the boundary. 
As in \cite[Section~3.1]{ERW22}, this is a moduli space
of smooth $N$-bundles with the fixed boundary identification
and a Spin structure on the vertical tangent bundle
extending the fixed boundary Spin structure.
Forgetting the relative Spin structure gives a weak
equivalence
\[
\rho:
\cC_{\Spin(16)}^2(\varnothing,S^{15})_N
\xrightarrow{\simeq}
B\Diff_\partial(N).
\]
Indeed, its homotopy fibre is the space of lifts
of the tangential classifying map $N\to BO(16)$ to
$B\Spin(16)$ extending the fixed boundary Spin structure, which is contractible since
$\partial N\to N$ is $2$-connected and
$B\Spin(16)\to BO(16)$ is $2$-coconnected. See \cite[Section~7.1]{ERW22} for details. 

Then we have the maps
\[
B\widehat{\Diff}_\partial(N)
\longrightarrow B\Diff_\partial(N)
\xleftarrow[\simeq]{\rho}
\cC_{\Spin(16)}^2(\varnothing,S^{15})_N
\]
which classify the underlying smooth bundles.
Fix a collared embedding of the Spin cobordism
$N:\varnothing\to S^{15}$. 
This choice determines basepoints
in the three moduli spaces in their embedding models. We denote the basepoint by $N$.

Finally, the tautological map
\[
\tau:
\cC_{\Spin(16)}^2(\varnothing,S^{15})_N
\longrightarrow
\Omega_{\varnothing,S^{15}}B\cC_{\Spin(16)}^2
\]
sends each cobordism to the path represented by its
$1$-simplex in the nerve of the cobordism category.
Here $\Omega_{\varnothing,S^{15}}B\cC_{\Spin(16)}^2$
denotes the space of paths in $B\cC_{\Spin(16)}^2$ from
the empty object to the sphere $S^{15}$ with its fixed
Spin structure. 

Choosing a based homotopy inverse $\rho^{-1}$, we define
$\widetilde\alpha_N$ by the composite
\begin{equation}\label{eq:punctured-PT}
	\widetilde\alpha_N: B\widehat{\Diff}_\partial(N)
	\rightarrow
	B\Diff_\partial(N)
	\xrightarrow{\rho^{-1}}
	\cC_{\Spin(16)}^2(\varnothing,S^{15})_N
	\xrightarrow{\tau}
	\Omega_{\varnothing,S^{15}}B\cC_{\Spin(16)}^2.
\end{equation}
Then $\widetilde\alpha_N$ sends the basepoint $N$ to the path
$\varnothing \rightsquigarrow S^{15}$ represented by the cobordism $N$.

Let $d_{\mathrm{cap}}:S^{15}\rightsquigarrow\varnothing$
be the path in $B\cC_{\Spin(16)}$ represented by
the chosen disk $D^{16}\subset\OP^2$ with induced Spin structure, and regarded as a cobordism
$S^{15}\to\varnothing$ using the fixed boundary identification
with $N$.
As in \cite[Section~7.3]{ERW22}, we have the following homotopy commutative diagram
\begin{equation}\label{eq:capping-square}
	\begin{tikzcd}[column sep=4em,row sep=large]
		B\widehat{\Diff}_\partial(N)
		\arrow[r,"B\operatorname{cap}"]
		\arrow[d,"i\circ\widetilde\alpha_N"'] &
		\Bdiff
		\arrow[d,"\operatorname{taut}"] \\
		\Omega_{\varnothing,S^{15}}\Bord
		\arrow[r,"{-*d_{\mathrm{cap}}}"'] &
		\Omega_\varnothing\Bord
	\end{tikzcd}
\end{equation}
Here the tautological map regards a closed $\OP^2$-bundle as
a family of cobordisms $\varnothing\to\varnothing$.
Under the Pontryagin--Thom equivalence,
$\operatorname{PT}\circ\operatorname{taut}\simeq\alpha_{\OP^2}$.
Thus $\alpha_N$ is represented by the composite
\[
\alpha_N: 
B\widehat{\Diff}_\partial(N)
\xrightarrow{B\operatorname{cap}}
\Bdiff
\xrightarrow{\operatorname{taut}}
\Omega_\varnothing\Bord
\xrightarrow[\simeq]{\operatorname{PT}}
\Omega^\infty MT\Spin(16).
\]
The homotopy commutative square \eqref{eq:capping-square} therefore gives
\begin{equation}\label{eq:homotopy-square}
	\alpha_N\simeq
	\operatorname{PT}\circ(-*d_{\mathrm{cap}})\circ i\circ\widetilde\alpha_N. 
\end{equation}
Normalizing at $N$, we have 
\[
\begin{aligned}
	\alpha_N(-)-\alpha_N(N)
	&\simeq
	\operatorname{PT}\bigl(i(\widetilde\alpha_N(-))*d_{\mathrm{cap}}\bigr)
	-\operatorname{PT}\bigl(i(\widetilde\alpha_N(N))*d_{\mathrm{cap}}\bigr)\\
	&\simeq \operatorname{PT}\Bigl(
	\bigl(i(\widetilde\alpha_N(-))*d_{\mathrm{cap}}\bigr)*
	\overline{\bigl(i(\widetilde\alpha_N(N))*d_{\mathrm{cap}}\bigr)}
	\Bigr)\\
	&\simeq
	\operatorname{PT}\Bigl(
	i(\widetilde\alpha_N(-))*\overline{i(\widetilde\alpha_N(N))}
	\Bigr),
\end{aligned}
\]
where the second homotopy is due to the fact that $\operatorname{PT}$ sends
loop concatenation to addition up to homotopy \cite{GMTW09},
and the last homotopy cancels
$d_{\mathrm{cap}}*\overline{d_{\mathrm{cap}}}$.

Consider the based homotopy equivalence
\begin{equation}\label{eq:Th}
	T=(-*\widetilde\alpha_N(N)):
	\Omega_\varnothing B\cC_{\Spin(16)}^2
	\xrightarrow{\simeq}
	\bigl(
	\Omega_{\varnothing,S^{15}}B\cC_{\Spin(16)}^2,
	\widetilde\alpha_N(N)
	\bigr),
\end{equation}
defined by concatenation with the path $\widetilde\alpha_N(N)$. Then 
\[
i(\widetilde\alpha_N(-))*\overline{i(\widetilde\alpha_N(N))}
\simeq
(\Omega i)\circ T^{-1}\circ\widetilde\alpha_N, 
\]
and we have
\begin{equation}\label{eq:cap-normalization}
	\alpha_N(-)-\alpha_N(N)
	\simeq \operatorname{PT}\circ(\Omega i)\circ T^{-1}
	\circ\widetilde\alpha_N.
\end{equation}
By \eqref{eq:z2},
$z_0=(\alpha_N(-)-\alpha_N(N))_*\xi_\partial$.
Applying \eqref{eq:cap-normalization} to
$\xi_\partial\otimes1$ gives
\[
z_0\otimes1=\operatorname{PT}_*(\Omega i)_*T_*^{-1}
(\widetilde\alpha_N)_*(\xi_\partial\otimes1).
\]
Combining this equality with \eqref{eq:spherical-image}
and $f=i\circ f_2$, and using that $\operatorname{PT}_*$
and $(\Omega i)_*$ are isomorphisms, gives

\begin{equation}\label{eq:Th-normalization}
	(\widetilde\alpha_N)_*(\xi_\partial\otimes1)
	=T_*(\Omega f_2)_*a_\Q.
\end{equation}

\subsection{Vanishing of the psc obstruction}

We now combine \eqref{eq:Th-normalization} with the
comparison square of Ebert--Randal-Williams to lift
a positive multiple of $\xi_\partial$ to the metric bundle.

\begin{proposition}\label{prop:integrallift}
There is a positive integer $k$ such that the class $k\xi_\partial$
lifts through
\[
 \pi_8\bigl(E\widehat{\Diff}_\partial(N)
 \times_{\widehat{\Diff}_\partial(N)}
 \cR^+(N)_{g_\partial}\bigr)
 \longrightarrow
 \pi_8 B\widehat{\Diff}_\partial(N).
\]
Consequently, $\delta_h(k\xi_\partial)=k\delta_h(\xi_\partial)=0$ and the family classified by $k\xi_\partial$ admits a continuous family of psc metrics with the fixed boundary condition. 
\end{proposition}
\begin{proof}

	Write
	\[
	\mathcal F_\partial:=\hofib_{S^{15}}(\Bpsc\xrightarrow{f_2} B\cC_{\Spin(16)}^2), 
	\]
	in which a point is given by a pair $(x,\gamma)$, where
	$x\in\Bpsc$ and $\gamma$ is a path from $f_2(x)$ to $S^{15}$
	in $B\cC_{\Spin(16)}^2$. In particular, there is a natural map
	\[
	\sigma:\Omega_{\varnothing,S^{15}}B\cC_{\Spin(16)}^2
	\longrightarrow\mathcal F_\partial,
	\qquad v\longmapsto(\varnothing,v).
	\]
	
	For the outgoing
	manifold $S^{15}$ and metric $g_\partial$, let $Y_\bullet$ be the semisimplicial space as in
	\cite[Definition~7.2.4]{ERW22}. 
	A point of $Y_p$ consists of levels
	$0<a_0<\cdots<a_p<1$, right stable psc metrics on
	$[a_{j-1},a_j]\times S^{15}$ for $1\leq j\leq p$,
	and an arbitrary psc metric on $[a_p,1]\times S^{15}$
	whose outgoing metric is $g_\partial$.
	These metrics have matching boundary values and
	product collars. 
	We have the following
	homotopy cartesian square, 
	\begin{equation}\label{eq:actionsquare}
		\begin{tikzcd}[column sep=2.5em,row sep=large]
			E\widehat{\Diff}_\partial(N)
			\times_{\widehat{\Diff}_\partial(N)}\cR^+(N)_{g_\partial}
			\arrow[r,"b"]\arrow[d]
			& \|Y_\bullet\|\arrow[d,"q"]\\
			B\widehat{\Diff}_\partial(N)
			\arrow[r,"\sigma\widetilde\alpha_N"']
			&\mathcal F_\partial
		\end{tikzcd}
	\end{equation}
	obtained from \cite[(7.2.1)]{ERW22}
	by restricting to cobordisms diffeomorphic to $N$
	relative to the boundary, and pulling back along
	$B\widehat{\Diff}_\partial(N)\to B\Diff_\partial(N)$.

	Evaluating the chosen homotopy in \eqref{eq:actionsquare}
	at $[e,h]$ gives the path
	\[
	\eta:q\circ b([e,h])
	\rightsquigarrow\sigma\widetilde\alpha_N(N).
	\]
	The lower map in \eqref{eq:actionsquare} sends
	$\xi_\partial$ to a torsion class. Indeed, consider the homotopy fibre sequence of $f_2$, followed by transport
	from $\varnothing$ to $S^{15}$ along $\widetilde\alpha_N(N)$, 
	\[
	\Omega\Bpsc\xrightarrow{\Omega f_2}
	\Omega B\cC_{\Spin(16)}^2
	\xrightarrow{\sigma T}\mathcal F_\partial.
	\]
	Applying $\sigma_*$ to \eqref{eq:Th-normalization} gives
	\begin{equation}\label{eq:psc-zero}
			\sigma_*(\widetilde\alpha_N)_*(\xi_\partial\otimes1)
			=\sigma_* T_*(\Omega f_2)_*a_\Q=0
	\end{equation}
	since $\sigma_* T_*(\Omega f_2)_*=0$ by exactness. Hence the integral class
	$(\sigma\widetilde\alpha_N)_*\xi_\partial$ is torsion, and we choose a positive $k$ that kills it. Concretely, let
	\[
	u:(S^8,*)\longrightarrow
	\bigl(B\widehat{\Diff}_\partial(N),N\bigr)
	\]
	represent $k\xi_\partial$, and choose a homotopy
	$H:S^8\times[0,1]\to\mathcal F_\partial$ satisfying
	\[
	\begin{aligned}
		H(s,0)=H(*,t)=\sigma\widetilde\alpha_N(N),\qquad
		H(s,1)=\sigma\widetilde\alpha_N(u(s)).
	\end{aligned}
	\]

	We first lift $u$ to the homotopy pullback
	$B\widehat{\Diff}_\partial(N)
	\times^h_{\mathcal F_\partial}\|Y_\bullet\|$, 
	whose points are triples $(x,y,\lambda)$,
	where $x\in B\widehat{\Diff}_\partial(N)$,
	$y\in\|Y_\bullet\|$, and
	$\lambda:[0,1]\to\mathcal F_\partial$ satisfies
	\[
	\lambda(0)=q(y),
	\qquad
	\lambda(1)=\sigma\widetilde\alpha_N(x).
	\]
	Indeed, the based lift is given by
	\[
	\begin{aligned}
		(S^8,*) &\longrightarrow
		\big(
		B\widehat{\Diff}_\partial(N)
		\times^h_{\mathcal F_\partial}\|Y_\bullet\|,
		(N,b([e,h]),\eta) 
		\big) \\
		s&\mapsto
		\bigl(u(s),b([e,h]),\eta*H(s,-)\bigr).
	\end{aligned}
	\]
	By the universal property, the homotopy in \eqref{eq:actionsquare} defines a weak equivalence
	\[
	E\widehat{\Diff}_\partial(N)
	\times_{\widehat{\Diff}_\partial(N)}
	\cR^+(N)_{g_\partial}
	\longrightarrow
	B\widehat{\Diff}_\partial(N)
	\times^h_{\mathcal F_\partial}\|Y_\bullet\|,
	\]
	Hence on $\pi_8$, the representative $u$ further lifts to the total space of the associated bundle:
	\begin{equation}\label{eq:u-tilde}
		\begin{aligned}
			\widetilde{u}: (S^8,*) \longrightarrow
			\big(
			E\widehat{\Diff}_\partial(N)
			\times_{\widehat{\Diff}_\partial(N)}
			\cR^+(N)_{g_\partial},
			[e,h] 
			\big). 
		\end{aligned}
	\end{equation}
	This proves the assertion.
\end{proof}

\subsection{Proof of Theorem A}

We cap the continuous metric family on the punctured bundle from
Proposition~\ref{prop:integrallift}, and smooth it to obtain the desired smooth fibrewise psc metric. 

\begin{proof}[Proof of
 Theorem~\ref{thm:A}]

 Let $k$ and the representative $u$ of $k\xi_\partial$
 be as in the proof of Proposition~\ref{prop:integrallift}.
 Let $E_\partial\to S^8$ be the smooth $N$-bundle
 classified by $u$ with its fixed boundary identification and vertical Spin structure. Let $\cR^+(E_\partial/S^8)_{g_\partial}$ be the pullback of \eqref{eq:metric-borel} along $u$, whose fibre at $x$ is
 $\cR^+((E_\partial)_x)_{g_\partial}$.
 \[
 \begin{tikzcd}[column sep=large,row sep=large]
 	\cR^+(E_\partial/S^8)_{g_\partial}
 	\arrow[r]\arrow[shift left, d]
 	&
 	E\widehat{\Diff}_\partial(N)
 	\times_{\widehat{\Diff}_\partial(N)}
 	\cR^+(N)_{g_\partial}
 	\arrow[d]
 	\\
 	S^8
 	\arrow[r,"u"']
 	\arrow[ur,dashed, "\widetilde{u}"]
 	\arrow[shift left, dashed, u,"m"]
 	&
 	B\widehat{\Diff}_\partial(N)
 \end{tikzcd}
 \]
 As in \eqref{eq:u-tilde}, the based lift $\widetilde u$ induces a section $m$ of the pullback bundle, and $m(x)$ is a psc metric on the $N$-torsor $(E_\partial)_x$ with the fixed boundary condition. In particular, $m(*)=h$.

Recall the capping homomorphism defined above
 \[
 \operatorname{cap}:\widehat{\Diff}_\partial(N)
 \to\widehat{\Diff}_0(\OP^2)
 \]
 The composition $(B\operatorname{cap})\circ u$ therefore classifies the smooth bundle
 \[
 E_+:=
 E_\partial\cup_{S^8\times S^{15}}
 (S^8\times D^{16})
 \xrightarrow{\pi}S^8,
 \]
which corresponds to the class $k\xi$, and
is isomorphic to a degree-$k$ pullback of the original Cayley plane family classified by $\xi$.

Gluing the chosen torpedo metric $g_{\mathrm{tor}}$ along the product collars gives a continuous map
\cite[Sections~2.1 and~2.3]{ERW22}
\begin{equation}\label{eq:cap}
	\begin{aligned}
		\operatorname{cap}_\mathrm{tor}:\cR^+(N)_{g_\partial}
		&\longrightarrow\cR^+(\OP^2),
		\\
		\mu&\longmapsto\mu\cup g_{\mathrm{tor}},
	\end{aligned}
\end{equation}
which is equivariant
with respect to the capping homomorphism $\operatorname{cap}$. 
Let $\cR^+(E_+/S^8)_{g_{\mathrm{tor}}}$ denote the bundle
of fibrewise psc metrics equal to $g_{\mathrm{tor}}$
on the fixed $D^{16}$-subbundle.
By the equivariance of $\operatorname{cap}_\mathrm{tor}$, we have a continuous section
\begin{equation}\label{eq:section-conti}
	S^8
	\xrightarrow{m}
	\cR^+(E_\partial/S^8)_{g_\partial}
	\xrightarrow{\operatorname{cap}_{\mathrm{tor}}}
	\cR^+(E_+/S^8)_{g_{\mathrm{tor}}}.
\end{equation}
Its value at the basepoint is
$h\cup g_{\mathrm{tor}}=g$.

Now it suffices for us to find a smooth replacement for the section \eqref{eq:section-conti}. 
Due to the psc condition, the typical fibre of
 $\cR^+(E_+/S^8)_{g_{\mathrm{tor}}}$
 is open in the affine Fréchet space
 \[
 g+
 \left\{
 a\in\Gamma^\infty(\OP^2;\Sym^2T^*\OP^2):
 a|_{D^{16}}=0
 \right\},
 \] 
 hence a locally convex manifold. 
 And since the transition diffeomorphisms of $E_+$
 fix $D^{16}$ and act smoothly on metrics by pullback,
 $\cR^+(E_+/S^8)_{g_{\mathrm{tor}}}$
 is a smooth locally trivial bundle with locally convex typical fibre.
 Then by \cite[Proposition~13]{Woc09}, the continuous section \eqref{eq:section-conti} is homotopic to a smooth one through sections fixing its value $g$
 at the basepoint. 
 
 Let $g_+$ be the smooth vertical Riemannian metric
 on $T_\pi E_+$ corresponding to this smooth section.
 Its restriction to each fibre has positive scalar curvature,
 and restriction to the trivial subbundle $S^8\times D^{16}$ coincides with $g_{\mathrm{tor}}$.
 By the choice of $g$ in Section~\ref{subsec:based-model},
$g_+|_{(E_+)_*}=g$ lies in 
$\cR^+_{\mathrm{std}}(\OP^2)$.
Since $S^8$ is path connected and $\Diff_0(\OP^2)$
preserves this component, the fibre metrics of $g_+$
define a section of the associated bundle with fibre
$\cR^+_{\mathrm{std}}(\OP^2)$.
The smooth family $(E_+,g_+)$ therefore gives a classifying map
\[
f_+:S^8\longrightarrow\mathcal M^+_{\mathrm{std}}.
\]
Evaluating $C$ on the underlying smooth bundle gives
\begin{equation}\label{eq:period}
	\begin{aligned}
		\langle (B\operatorname{cap}\circ u)^*C,[S^8]\rangle=k\langle C,\xi\rangle\ne0.
	\end{aligned}
\end{equation}

The fibrewise psc metric $g_+$ makes the vertical
Dirac operators invertible, hence
$\pi_!^{KO}(1)=0\in KO^{-16}(S^8)$.
By composition of Spin pushforwards,
\[
\alpha_{KO}(E_+)
=(S^8\to *)_!^{KO}\circ\pi_!^{KO}(1)=0.
\]
By connectivity, $E_+$ is a $24$-dimensional String manifold. 
Then $\Wit(E_+)$ is a modular form of weight $12$ with constant term $\Ahat(E_+)=0$, hence a scalar multiple of $\Delta$.
Since $TS^8$ is stably trivial, \eqref{eq:total} gives
\[
[q]\Wit(E_+)=
k\langle C,\xi\rangle.
\]
As $[q]\Delta=1$, we conclude that
\[
\Wit(E_+)=k\langle C,\xi\rangle\Delta\ne0.
\]
By \eqref{eq:homogeneouskernel}, the bundle $E_+$
is not homogeneous.
\end{proof}


\begin{thebibliography}{HHLZ22}
	\setlength{\itemsep}{2pt plus 1pt}
	
	\bibitem[AS71V]{AS71V}
	M. F. Atiyah and I. M. Singer,
	\emph{The index of elliptic operators. V},
	Ann. of Math. (2) \textbf{93} (1971), 139--149.
	\url{https://doi.org/10.2307/1970757}.
	
	\bibitem[BH58]{BH58}
	A. Borel and F. Hirzebruch,
	\emph{Characteristic classes and homogeneous spaces, I},
	Amer. J. Math. \textbf{80} (1958), 458--538.
	\url{https://doi.org/10.2307/2372795}.
	
	\bibitem[BHSW10]{BHSW10}
	B. Botvinnik, B. Hanke, T. Schick, and M. Walsh,
	\emph{Homotopy groups of the moduli space of metrics of positive scalar curvature},
	Geom. Topol. \textbf{14} (2010), 2047--2076.
	\url{https://doi.org/10.2140/gt.2010.14.2047}.
	
	\bibitem[Che04]{Che04}
	V. Chernysh,
	\emph{On the homotopy type of the space $\mathcal{R}^+(M)$},
	arXiv:math/0405235 (2004).
	\url{https://arxiv.org/abs/math/0405235}.
	
	\bibitem[Des09]{Des09}
	A. Dessai,
	\emph{Some geometric properties of the Witten genus},
	in \emph{Alpine Perspectives on Algebraic Topology},
	Contemp. Math., vol.~504, Amer. Math. Soc., Providence, RI, 2009,
	pp.~99--115.
	\url{https://doi.org/10.1090/conm/504/09877}.
	
	\bibitem[ERW22]{ERW22}
	J. Ebert and O. Randal-Williams,
	\emph{The positive scalar curvature cobordism category},
	Duke Math. J. \textbf{171} (2022), 2275--2406.
	\url{https://doi.org/10.1215/00127094-2022-0023}.
	
	\bibitem[EW24]{EW24}
	J. Ebert and M. Wiemeler,
	\emph{On the homotopy type of the space of metrics of positive scalar curvature},
	J. Eur. Math. Soc. \textbf{26} (2024), 3327--3363.
	\url{https://doi.org/10.4171/JEMS/1333}.
	
	\bibitem[Fre25]{Fre25}
	G. Frenck,
	\emph{Characteristic numbers of manifold bundles over spheres and positive curvature via block bundles},
	with an appendix by G. Frenck and J. Reinhold,
	Forum Math. Sigma \textbf{13} (2025), e138.
	\url{https://doi.org/10.1017/fms.2025.10054}.
	
	\bibitem[GMTW09]{GMTW09}
	S. Galatius, I. Madsen, U. Tillmann, and M. Weiss,
	\emph{The homotopy type of the cobordism category},
	Acta Math. \textbf{202} (2009), 195--239.
	\url{https://doi.org/10.1007/s11511-009-0036-9}.
	
	\bibitem[GRW14]{GRW14MMM}
	S. Galatius and O. Randal-Williams,
	\emph{Detecting and realising characteristic classes of manifold bundles},
	in \emph{Algebraic Topology: Applications and New Directions},
	Contemp. Math., vol.~620, Amer. Math. Soc., Providence, RI, 2014,
	pp.~99--110.
	\url{https://doi.org/10.1090/conm/620/12365}.
	
	\bibitem[GRW17]{GRW17}
	S. Galatius and O. Randal-Williams,
	\emph{Homological stability for moduli spaces of high dimensional manifolds. II},
	Ann. of Math. (2) \textbf{186} (2017), 127--204.
	\url{https://doi.org/10.4007/annals.2017.186.1.4}.
	
	\bibitem[GRW20]{GRWGuide}
	S. Galatius and O. Randal-Williams,
	\emph{Moduli spaces of manifolds: a user's guide},
	in \emph{Handbook of Homotopy Theory} (H. Miller, ed.),
	Chapman \& Hall/CRC, Boca Raton, FL, 2020, pp.~443--485.
	\url{https://doi.org/10.1201/9781351251624-12}.
	
	\bibitem[HH22]{HH22}
	F. Han and R. Huang,
	\emph{On characteristic numbers of 24 dimensional string manifolds},
	Math. Z. \textbf{300} (2022), 2309--2331.
	\url{https://doi.org/10.1007/s00209-021-02877-6}.
	
	\bibitem[HHLZ22]{HHLZ21}
	F. Han, R. Huang, K. Liu, and W. Zhang,
	\emph{Cubic forms, anomaly cancellation and modularity},
	Adv. Math. \textbf{394} (2022), 108023.
	\url{https://doi.org/10.1016/j.aim.2021.108023}.
	
	\bibitem[Hit74]{Hit74}
	N. Hitchin,
	\emph{Harmonic spinors},
	Adv. Math. \textbf{14} (1974), 1--55.
	\url{https://doi.org/10.1016/0001-8708(74)90021-8}.
	
	\bibitem[KS93]{KS93}
	M. Kreck and S. Stolz,
	\emph{$\mathrm{HP}^2$-bundles and elliptic homology},
	Acta Math. \textbf{171} (1993), 231--261.
	\url{https://doi.org/10.1007/BF02392533}.
	
	\bibitem[MS26]{MS26}
	S. Mandal and A. Sarkar,
	\emph{Relative smooth surgery structure sets of thickenings of the Cayley projective plane and applications},
	arXiv:2607.20362 (2026).
	\url{https://arxiv.org/abs/2607.20362v2}.
	
	\bibitem[McT14]{McT14}
	C. McTague,
	\emph{The Cayley plane and string bordism},
	Geom. Topol. \textbf{18} (2014), 2045--2078.
	\url{https://doi.org/10.2140/gt.2014.18.2045}.
	
	\bibitem[MM65]{MM65}
	J. W. Milnor and J. C. Moore,
	\emph{On the structure of Hopf algebras},
	Ann. of Math. (2) \textbf{81} (1965), 211--264.
	\url{https://doi.org/10.2307/1970615}.
	
	\bibitem[St92]{St92}
	S. Stolz,
	\emph{Simply connected manifolds of positive scalar curvature},
	Ann. of Math. (2) \textbf{136} (1992), 511--540.
	\url{https://doi.org/10.2307/2946598}.
	
	\bibitem[Stong68]{Stong68}
	R. E. Stong,
	\emph{Notes on Cobordism Theory},
	Mathematical Notes, no.~7, Princeton University Press, Princeton, NJ, 1968.
	
	\bibitem[Woc09]{Woc09}
	C. Wockel,
	\emph{A generalization of Steenrod's approximation theorem},
	Arch. Math. (Brno) \textbf{45} (2009), 95--104.
	\url{https://dml.cz/handle/10338.dmlcz/128295}.
	
\end{thebibliography}
\end{document}